\documentclass{amsart}
\usepackage{hyperref}

\usepackage{amsrefs,amsmath,amsfonts,amssymb,graphicx,latexsym,color,amsthm,euscript}

\newtheorem{theorem}{Theorem}
\newtheorem{proposition}[theorem]{Proposition}

\newtheorem{lemma}[theorem]{Lemma}

\newcommand{\eqdef}{\overset{\mbox{\tiny{def}}}{=}}
\newcommand{\dis}{\displaystyle}

\newcommand{\rmi}{{\rm i}}
\newcommand{\rmre}{{\rm Re}}

\newcommand{\ang}[1]{ \left< {#1} \right> }

\newcommand{\R}{\mathbb{R}}

\newcommand{\fcnP}{\hat{\rho}}

\newcommand{\ind}{ {\mathbf 1}}

\newcommand{\FP}{\mathbf{P}}
\newcommand{\FL}{L}
\newcommand{\FI}{\mathbf{I}}

\newcommand{\CE}{\mathcal{E}}
\newcommand{\CF}{\mathcal{F}}

\newcommand{\Ndim}{3}
\newcommand{\threed}{{\mathbb R}^\Ndim}

\newcommand{\na}{\nabla}
\newcommand{\la}{\lambda}
\newcommand{\de}{\delta}
\newcommand{\pa}{\partial}

\newcommand{\Ga}{\Gamma}

\newcommand{\highG}{\Theta}
\newcommand{\highB}{\Lambda}

\newcommand{\sourceG}{g}
\newcommand{\solU}{f}
\newcommand{\vel}{p}

\newcommand{\wN}{\ell}

\newcommand{\wE}{j}
\newcommand{\wK}{\delta}

\numberwithin{equation}{section}
\numberwithin{theorem}{section}

\begin{document}

\title[Decay of the Soft Potential relativistic Boltzmann equation in $\R^3_x$]
{Large-Time Decay of the Soft Potential relativistic 
Boltzmann equation in $\R^3_x$}

\author[R. M. Strain]{Robert M. Strain}
\thanks{R.M.S. was partially supported by the NSF grant DMS-0901463, and an Alfred P. Sloan Foundation Research Fellowship.}

\author[K. Zhu]{Keya Zhu}
\address{University of Pennsylvania, Department of Mathematics, David Rittenhouse Lab, 209 South 33rd Street, Philadelphia, PA 19104-6395, USA} 
\email{strain at math.upenn.edu  \& zhuk at math.upenn.edu}
\urladdr{http://www.math.upenn.edu/~strain/  \&   http://www.math.upenn.edu/~zhuk/}

\begin{abstract}
For the relativistic Boltzmann equation in $\R^3_x$, this work proves the global existence, uniqueness, positivity, and optimal time convergence rates to the relativistic Maxwellian for solutions which start out sufficiently close under the general physical soft potential assumption proposed in 1988 \cite{MR933458}.
\end{abstract}


\setcounter{tocdepth}{1}

\maketitle

\thispagestyle{empty}

\section{Introduction and statement of the main results}

In the study of fast moving particles the relativistic Boltzmann equation is a fundamental physical model \cite{MR1898707,MR635279,MR1379589}.  In early work of Glassey \& Strauss \cite{MR1105532,MR1211782,MR1321370} from 1991, 1993, and 1995 global existence and uniqueness of nearby equilibrium solutions, and their large time convergence rates were shown in the torus ($\mathbb{T}^3_x$) and also in the whole space ($\mathbb{R}^3_x$).  On the torus the convergence rates are exponentially fast, and in the whole space the convergence rates are polynomial.  Their assumptions on the differential cross-section, $\sigma(g,\theta)$, fell into the regime of hard potentials.  Further results for the hard potential case can be found in \cite{MR2249574}.  However, for relativistic interactions, when one considers particles that are fast moving, a very important physical regime is the soft potential case; see \cite{DEnotMSI} for a physical point of view.  In recent work \cite{Strain2010}, the global existence, uniqueness, and rapid time convergence rates for nearby equilibrium solutions to the soft potential relativistic Boltzmann equation was shown on the torus ($\mathbb{T}^3_x$).  However the difficult whole space case has remained a challenging open problem.  In this work we prove the global existence, uniqueness, and optimal time convergence rates for nearby equilibrium solutions to the relativistic Boltzmann equation in $\R^3_x$ under the general physical soft potential assumption proposed in 1988 by \cite{MR933458}.

The relativistic Boltzmann equation is given by
\begin{equation}
\partial _{t}F+\hat{p}\cdot \nabla_x F=\mathcal{Q}(F,F).
\label{RBF}
\end{equation}
The solution, $F=F(t,x,p)$, is a  function of time $t\in [0,\infty)$, space $x\in \mathbb{R}^3$ and momentum $p\in \mathbb{R}^3$.  
It is conventional to denote the normalized velocity as
\begin{equation}
\hat{p}=\frac{p}{p^0}=\frac{p}{\sqrt{1+|p|^2}}.
\label{normV}
\end{equation}
Steady states of the relativistic Boltzmann equation are the  J{\"{u}}ttner solutions, which are commonly called relativistic Maxwellians.
They are given by
 \begin{equation}
J (p)
=
\frac{e^{- p^0}}{4\pi}.
\label{juttner}
\end{equation}
The collision operator, $\mathcal{Q}(F,F)$, is defined in \eqref{collisionCM}.  For the sake of simplicity but without loss of generality we have normalized the physical constants to be one.

The entropy of the relativistic Boltzmann equation is physically defined as
\begin{equation}
\mathcal{H}(t)
\eqdef
-\int_{\mathbb{R}^3}  dx ~
\int_{\mathbb{R}^3}   dp ~
 F(t, x, p)\ln F(t, x, p).  
\nonumber
\end{equation}
Boltzmann's H-Theorem then corresponds to the formal differential inequality
\begin{equation}
\frac d{dt}\mathcal{H}(t)
 \ge 0,  \nonumber
\end{equation}
which predicts that the entropy of solutions will be non-decreasing as time passes.   It is well known that the steady state relativistic Maxwellians \eqref{juttner} maximize the entropy which
grants the intuition of convergence to \eqref{juttner} in large time.

It is this physical reasoning that our main results, as stated below, make mathematically rigorous in the context of   perturbations of the relativistic Maxwellian for a general class of soft potential cross-sections in the whole space 
$\mathbb{R}^3_x$.

\subsection{Notation}
In this section we define several notations which will be used throughout the article.  We first introduce the center of momentum expression for the collision operator, as presented in \cite{sKRM2011}.  
In particular we have
\begin{equation}
\mathcal{Q}(f,h)
=
\int_{\mathbb{R}^3}  dq ~
\int_{\mathbb{S}^{2}} d\omega ~
~ v_{\o} ~ \sigma (g,\theta ) ~ [f(p^{\prime })h(q^{\prime})-f(p)h(q)].
\label{collisionCM}
\end{equation}
where $v_{\o}=v_{\o}(p,q)$ is the M{\o}ller velocity given by
\begin{equation}
v_{\o}=
v_{\o}(p,q)
\eqdef
\sqrt{\left| \frac{p}{p^0}-\frac{q}{q^0}\right|^2-\left| \frac{p}{p^0}\times\frac{q}{q^0}\right|^2}
=
\frac{ g\sqrt{s}}{p^0 q^0}.
\label{moller}
\end{equation}
Here a relativistic particle has momentum $p=(p^1, p^2 , p^\Ndim)\in\threed$, with 
its energy defined by $p^{0}=\sqrt{1+|p|^2}$ where $|p|^2 \eqdef p\cdot p$.   The post-collisional momentum in the expression (\ref{collisionCM}) can then be written as:
\begin{equation}
\begin{split}
p^\prime&=\frac{p+q}{2}+\frac{g}{2}\left(\omega+(\gamma-1)(p+q)\frac{(p+q)\cdot \omega}{|p+q|^2}\right), 
\\
q^\prime&=\frac{p+q}{2}-\frac{g}{2}\left(\omega+(\gamma-1)(p+q)\frac{(p+q)\cdot \omega}{|p+q|^2}\right),
\label{postCOLLvelCMsec2}
\end{split}
\end{equation}
where $\gamma =(p^0+q^0)/\sqrt{s}$.
These will satisfy \eqref{collisionalCONSERVATION}.  The angle further satisfies
$
\cos\theta = k\cdot \omega
$
with $k = k(p,q)$ and $| k| = 1$.    The unit vector, $k$, has a complicated expression 
as given in \cite[Eq. (14)]{sKRM2011} but its precise form will be inessential.

Now conservation of momentum and energy is given as
\begin{equation}
\begin{split}
p^{0}+q^{0} & =p^{\prime 0}+q^{\prime 0},
\\
p+q & =p^{\prime}+q^{\prime}.
\end{split}
\label{collisionalCONSERVATION}
\end{equation}
Furthermore, the relative momentum, $g$,  is denoted
\begin{gather}
g 
\eqdef
\sqrt{2(p^0 q^0 - p\cdot q-1)}.
\label{gDEFINITION}
\end{gather}
Then the quantity ``$s$'' is defined as
\begin{eqnarray}
s 
\eqdef
 2(p^0 q^0 - p\cdot q+1)\ge 0.
\label{sDEFINITION}
\end{eqnarray}
Notice that $s=g^2+4$.    Also $\sigma(g, \theta)$ is the differential cross-section or scattering kernel; it is designed to measure the interactions between particles.     We give a standard warning to the reader that this notation, which is used in \cite{MR635279}, sometimes differs from other authors notation by a constant factor.

For an integrable function $g: \R^3_x\to\R$, its Fourier transform
is defined by
\begin{equation*}
  \widehat{g}(k)= \CF g(k) \eqdef \int_{\R^3} e^{-2\pi i x\cdot k} g(x)dx, \quad
  x\cdot
   k\eqdef\sum_{j=1}^3 x_jk_j,
   \quad
   k\in\R^3,
\end{equation*}
where $i =\sqrt{-1}\in \mathbb{C}$ is the imaginary unit. For two
complex vectors $a,b\in\mathbb{C}^3$, $(a\mid b)=a\cdot
\overline{b}$ denotes the dot product over the complex field, where
$\overline{b}$ is the ordinary complex conjugate of $b$.

For a function $f(t,x,p)$ with $x\in\mathbb{R}^3_x$, $p\in\mathbb{R}^3_p$, $t\in [0,\infty)$, 
we define the norm
$$
\|f\|_{L^{r_1}_tL^{r_2}_pL^{r_3}_x}
\eqdef
\left(\int^\infty_0\left(\int_{\mathbb{R}^3}\left(\int_{\mathbb{R}^3}|f(t,x,p)|^{r_3}\,dx\right)^{\frac{r_2}{r_3}}\,dp\right)^{\frac{r_1}{r_2}}\,dt\right)^{\frac{1}{r_1}}.
$$
Here $r_i \in [1, \infty)$.  
In the above expression, it is understood that if $r_i=\infty$ for some $i$, then the expression is modified accordingly by replacing the $L^{r_i}$ norm by the $L^\infty$ norm as usual.  In particular 
for any $r_i \in [1, \infty]$
we have 
$$
\|f \|_{L^{r_1}_tL^{r_2}_pL^{r_3}_x}
=
\left\|  \left\|  \|f \|_{L^{r_3}(\R^3_x)} \right\|_{L^{r_2}(\R_p^3)}   \right\|_{L^{r_1}([0,\infty))}.
$$
 The norm $\| f \|_{L^{r_2}_pL^{r_3}_x}$ is defined similarly (and it is a function in $t$). 
 The norm $\| f \|_{L^{r_3}_x}$ is also defined similarly (and it is analogously function in $t$ and $p$). 
 
To study the linear and non-linear time decay rates we introduce $\sigma_{r,m}$ as
 \begin{equation}
\label{rateLIN}
    \sigma_{r,m}\eqdef \frac{\Ndim}{2}\left(\frac{1}{r}-\frac{1}{2}\right)+\frac{m}{2}, 
    \quad 
    m \ge 0.
\end{equation}
Above the parameter $r$ satisfies $r\in [1,2]$.  This $\sigma_{r,m}$ is the standard notation for studying time decay rates for kinetic equations in the whole space.
 
 We also use the norms
$\|\cdot\|_{L^2_\vel(\dot{H}^m_x)}$ and $\|\cdot\|_{L^2_\vel(H^m_x)}$ with $m\geq 0$.
In this definition $\dot{H}^m_x=\dot{H}^m(\threed_x)$ is the standard homogeneous $L^2_x$ based Sobolev space. 
The $L^2(\mathbb{R}^3_p)$ inner product in the momentum variable $p$ is denoted $\langle \cdot, \cdot \rangle$.
Now, for $\ell\in\mathbb{R}$, we define the following weight function
\begin{equation}
w_{\ell}= w_\ell(p)
\eqdef 
(p^0)^{\ell b/2}.
\label{weight}
\end{equation}
The constant $b>0$ is defined for the soft-potentials in \eqref{hypSOFT} below.  
For the soft potentials, we will observe later on that  $w_1(p) \approx 1/\nu(p)$ (Lemma \ref{nuEST}). 

We also define the important temporal weight:   
\begin{equation}
\varpi_k
=
\varpi_k(t)
\eqdef
(1+t)^k, \quad k \ge 0.
\label{TIMEweight}
\end{equation}
Lastly, the notation $A \lesssim B$ will imply that a positive constant $C$ exists such that $A \leq C B$ holds uniformly over the range of parameters which are present in the inequality and  moreover that the precise magnitude of the constant is unimportant.  The notation $B \gtrsim A$ is equivalent to $A \lesssim B$, and $A \approx B$ means that both inequalities $A \lesssim B$ and $B \lesssim A$ hold simultaneously.   

Throughout this paper, furthermore $C$  denotes
some positive (generally large) constant and $c$ denotes some positive (generally small) constant, where both $C$ and
$c$ may take different values in different places.


\subsection{Perturbation equation and statement of the main results}

We will now explain in detail the main results of this paper.  
We define the standard perturbation $f(t,x,p)$ to the relativistic Maxwellian \eqref{juttner} as 
$$
F \eqdef J +\sqrt{J } f. 
$$
With \eqref{collisionalCONSERVATION} we observe that the quadratic collision operator \eqref{collisionCM} satisfies
$$
\mathcal{Q}(J, J)
=0.
$$
 Then the relativistic Boltzmann equation \eqref{RBF} which will be satisfied by the perturbation $f=f(t,x,p)$ is given by
\begin{gather}
 \partial_t f + \hat{p}\cdot \nabla_x f + L (f)
=
\Gamma (f,f),
\quad
f(0, x, p)=f_0(x,p).
\label{rBoltz0}
\end{gather}
 The linear operator $L( f)$ is defined in (\ref{L}).  And the non-linear operator $\Gamma(f,f)$ is defined  in (\ref{gamma0}).  They are derived from an expansion of the relativistic Boltzmann collision operator \eqref{collisionCM}.
In particular, the linearized collision operator is given by 
\begin{gather}
 L(h)
 \eqdef 
 -J^{-1/2}\mathcal{Q}(J ,\sqrt{J} h)- J^{-1/2}\mathcal{Q}(\sqrt{J} h,J)
 \label{L}
 \\
  =
 \nu(p) h-K (h).
 \notag
\end{gather}
Above the multiplication operator takes the form
\begin{equation}
\nu(p) \eqdef
\int_{\mathbb{R}^3} ~  dq 
\int_{\mathbb{S}^{2}} ~ d\omega
~ v_{\o} ~ \sigma(g,\theta)~ J(q).
\label{nuDEF}
\end{equation}
The remaining integral operator is
\begin{gather}
K(h) 
\eqdef
\int_{\mathbb{R}^3} ~  dq 
\int_{\mathbb{S}^{2}} ~ d\omega
~ v_{\o} ~ \sigma(g,\theta)
\sqrt{J(q)}\left\{ \sqrt{J(q^{\prime})} ~ h(p^{\prime })+\sqrt{J(p^{\prime})} ~ h(q^{\prime})\right\}
\notag
\\
-\int_{\mathbb{R}^3} ~  dq 
\int_{\mathbb{S}^{2}} ~ d\omega
~ v_{\o} ~ \sigma(g,\theta)
~ \sqrt{J(q) J(p)} ~ h(q)
\label{compactK}
\\
=
K_2(h) - K_1(h).
\notag
\end{gather}
The non-linear part of the collision operator is defined as
\begin{multline}
\Gamma (h_1,h_2)
\eqdef
 J^{-1/2}\mathcal{Q}(\sqrt{J} h_1,\sqrt{J} h_2)
\label{gamma0}
\\
=
\int_{\mathbb{R}^3} ~  dq 
\int_{\mathbb{S}^{2}} ~ d\omega
~ v_{\o} ~ \sigma(g,\theta)~
 \sqrt{J(q)}
 [h_1(p^{\prime })h_2(q^{\prime})-h_1(p)h_2(q)].
\end{multline}
For these operators, we use the following general conditions on the kernel. 

\subsection*{\bf Hypothesis on the collision kernel:}   
{\it For soft potentials we assume the collision kernel in \eqref{collisionCM} satisfies 
the following growth/decay estimates 
\begin{equation}
\begin{split}
\sigma (g,\theta) 
 &
\lesssim  g^{-b} ~ \sigma_0(\theta),
\\ 
\sigma (g,\theta) 
 &
\gtrsim
\left( \frac{g}{\sqrt{s}}\right) g^{-b} ~ \sigma_0(\theta).
\end{split}
\label{hypSOFT}
\end{equation}
We consider angular factors
$0 \le \sigma_0(\theta) \lesssim \sin^\gamma\theta$
with $\gamma> -2$.  Additionally  $\sigma_0(\theta)$ should be non-zero on a set of positive measure.
We suppose $0 < b <\min(4,4+\gamma)$.

For hard potentials we make the assumption
\begin{equation}
\sigma (g,\theta) 
 \lesssim
  \left(  g^{a}+ g^{-b}\right) ~ \sigma_0(\theta) ,
\quad
\sigma (g,\theta) 
 \gtrsim
  \left( \frac{g}{\sqrt{s}}\right) g^{a}  ~ \sigma_0(\theta).
\label{hypHARD}
\end{equation}
In addition to the previous parameter ranges we consider $0\le a\le 2+\gamma$ and also $0 \le b <\min(4,4+\gamma)$ (in this case we allow the possibility of $b=0$).
}

\bigskip

This hypothesis contains the general physical assumption on the kernel which was introduced in \cite{MR933458}  (and we add the corresponding necessary lower bounds).
 
We are now ready to state our main result:

\begin{theorem}
\label{mainGlobal}
Choose 
$\ell\geq \max\{0, 3/b-1\}$, 
$r\in [1,6/5)$ and $k\in (1/2, \sigma_{r,0}]$ where $\sigma_{r,0}$ is given by \eqref{rateLIN}.  Consider initial data $f_0 = f_0(x,p) \in L^2_p L^r_x \cap L^\infty_p(L^2_x\cap L^\infty_x)$.
  There is an $\eta>0$ such that if $\|f_0\|_{L^2_p L^r_x}+\|w_{\ell+k}f_0\|_{L^\infty_p(L^2_x\cap L^\infty_x)} \le \eta$, then   there exists a unique global in time mild solution \eqref{nonLinProbSol}, $f = f(t,x,p)$, to the relativistic Boltzmann equation \eqref{rBoltz0} with a soft potential kernel \eqref{hypSOFT} which satisfies
$$
\|w_\ell f\|_{L^\infty_p(L^2_x\cap L^\infty_x)}(t) \le C_{\ell, k}(1+t)^{-k} \left(\|f_0\|_{L^2_p L^r_x}+\|w_{\ell+k}f_0\|_{L^\infty_p(L^2_x\cap L^\infty_x)}\right).
$$ 
These solutions are continuous if it is so initially.  We furthermore have the positivity, in other words $F= \mu + \sqrt{\mu} f \ge 0$, if $F_0= \mu + \sqrt{\mu} f_0 \ge 0$.
\end{theorem}

In the next sub-section we will discuss some historical results related to our main theorem to explain what has been done in the past in connection with our result.  Explanations about why the convergence rates in Theorem \ref{mainGlobal} are said to be optimal can be found for instance in \cite{DuanStrainCPAM2011,sNonCutOp}.

\subsection{Historical discussion}
The relativistic Boltzmann equation is the primary model in relativistic collisional Kinetic theory.  In the next few paragraphs, we will provide a short review of the mathematical theory of this equation.  We mention a few books on relativistic Kinetic theory as for instance \cite{MR1898707,MR635279,MR1379589}.

In 1988, Dudy{\'n}ski and Ekiel-Je{\.z}ewska \cite{MR933458} proved that the linear relativistic Boltzmann equation admits unique solutions in $L^2$.  Afterwards, Dudy{\'n}ski \cite{MR1031410} studied the long time and small-mean-free-path limits of these solutions.

In the context of large data global in time weak solutions, the theory of DiPerna-Lions \cite{MR1014927} renormalized solutions was extended  to the relativistic Boltzmann equation in 1992 also by Dudy{\'n}ski and Ekiel-Je{\.z}ewska \cite{MR1151987}.  This result uses the causality of the relativistic Boltzmann equation 
\cite{MR818441,MR841735}.  
Results on the regularity of the gain term are given in \cite{MR1402446,MR1480243};
 the strong $L^1$ compactness is studied by  Andr{\'e}asson \cite{MR1402446}.
These are generalizations of Lions \cite{MR1284432} result in the non-relativistic case. 
Further developments on renormalized weak solutions can be found in 
\cite{MR1714446,MR1676150}.

Notice also the studies of the Newtonian limit \cite{MR2098116,MR2679588} for the Boltzmann equation.  
We further mention theories of unique global in time solutions with initial data that is near Vacuum as in 
Glassey \cite{MR2217287}
and 
\cite{MR2217287,MR2679588}. 

Note further the study of the collision map
 and the pre-post collisional change of variables from \cite{MR1105532}.  
Then \cite{MR2543323} provides uniform $L^2$-stability estimates for the relativistic
   Boltzmann equation.  
 Now there is a mathematically rigorous result connecting the relativistic Euler equations to the relativistic Boltzmann equation via the Hilbert expansion as in \cite{ssHilbert}.

We point out results on global existence of unique smooth solutions which are initially close to the relativistic Maxwellian for the 
relativistic Landau-Maxwell system \cite{MR2100057}, and then for the relativistic Landau \cite{MR2289548} equation as well.  Further \cite{MR2514726} proves the smoothing effects for relativistic Landau-Maxwell system.  
And \cite{MR2593052} proves time decay rates in whole space for the relativistic Boltzmann equation (with certain hard potentials) and the relativistic Landau equation as well.

Further previous results for unique strong solutions to the hard potential relativistic Boltzmann equation are as follows.  In 1993 
Glassey and Strauss \cite{MR1211782} proved asymptotic stability similar to  Theorem \ref{mainGlobal}  in  
$L^\infty_\ell$ with $\ell>3/2$ in $\mathbb{T}^3_x$.  They
consider collisional cross-sections which satisfy \eqref{hypHARD} for the parameters
$b\in [0,1/2)$, $a\in [0, 2-2b)$ and either $\gamma \ge 0$ or 
$$
\left| \gamma \right| 
< 
\min\left\{2-a, \frac{1}{2} -b, \frac{2-2b - a}{3}
\right\},
$$
which restricts to $\gamma > -\frac{1}{2}$ if $b=0$ say.
They further assume a related growth bound on the derivative of the  cross-section
$
\left| \frac{\partial \sigma}{\partial g} \right|.
$
In \cite{MR2249574} this growth bound was removed while the rest of the assumptions on the cross-section from \cite{MR1211782} remained the same.   These results also sometimes work in smoother function spaces, and we note that we could also include space-time regularity to our solutions spaces.   For the Cauchy problem, under similar assumptions on the collisional cross section,
Glassey \& Strauss \cite{MR1321370} in 1995
proved the global existence and uniqueness of nearby equilibrium solutions, and their large time polynomial convergence rates.  

However, it has been noted that for relativistic interactions, when one considers particles that are fast moving, a very important physical regime is the soft potential case \cite{DEnotMSI}.  Recently \cite{Strain2010}, the global existence, uniqueness, and rapid time convergence rates for nearby equilibrium solutions to the soft potential relativistic Boltzmann equation was shown on the torus ($\mathbb{T}^3_x$).  However the difficult problem of determining the optimal convergence rates in the whole space case under the full soft potential assumption from \eqref{hypSOFT} remained an open problem open prior to the main theorem of this paper.    We also believe that the methods used in this paper can be used to to treat the full hard potential assumption from \eqref{hypHARD}, however we consider it to be worthwhile to carry out this extension.

We reference other related and important previous work such as 
\cite{MR2013332,MR1211782,MR1105532,MR2209761,MR2366140,Mouhot:1173020,gsNonCutJAMS,gsNonCutADVMATH,gsNonCutA,MR2679358,MR0259662,MR2284213,MR2082240,
MR2116276,MR576265,MR575897,MR677262,MR1379589,MR1361017,Strain2010,sNonCutOp}.  We refer the reader to the discussions in \cite{sNonCutOp,DuanStrainCPAM2011,MR2754344} for more detailed explanations of historical developments and mathematical methods used.  We will discuss key new ideas and approaches during the course of the paper at the appropriate mathematical places below.

The remainder of this article is organized as follows.
In Section \ref{sec:22decay} we give the proof of  $L^2_{p}L^2_x$ decay for solutions to the linearized equation \eqref{rBoltz}.
Then in Section \ref{I2decay} we prove the linear $L^\infty_{p}L^2_x$ time decay.
Section \ref{sec:IIdecay} briefly explains the linear $L^\infty_{p}L^\infty_x$ time decay.
Lastly, in Section \ref{sec:NLdecay} we use the results from the previous sections to establish the non-linear time decay rates and the global existence of solutions.

\section{Linear decay theory in $L^2_{p}L^2_x$}\label{sec:22decay}

Our main goal in this section is to prove time decay of solutions to the linearized relativistic Boltzmann equation with the soft potentials \eqref{hypSOFT} in $L^2_{p}L^2_x$.  We consider the linearized equation with a microscopic source $\sourceG=\sourceG(t,x,\vel)$:
\begin{equation}\label{ls}
    \left\{\begin{array}{l}
  \dis     \pa_t \solU+\hat{p}\cdot\na_x \solU +\FL \solU =\sourceG,\\
\dis \solU|_{t=0}=\solU_0.
    \end{array}\right.
\end{equation}
For the nonlinear system \eqref{rBoltz0}, the non-homogeneous source takes the form of
\begin{equation}\label{def.g.non}
    \sourceG=\Ga(\solU,\solU).
\end{equation}
In this case $\sourceG=\{\FI-\FP\} \sourceG$, is microscopic as in \eqref{hydro}.  Then it is standard to observe that solutions of \eqref{ls} formally take the following form
\begin{equation}
    \solU(t)=U(t)\solU_0 + \int_0^t ds ~U(t-s)~\sourceG(s), \quad 
    U(t)\eqdef e^{-t\left( \FL+\vel\cdot\na_x \right)}.
    \label{ls.semi}
\end{equation}
Here $U(t)$ is the linear solution operator for the Cauchy problem corresponding to  \eqref{ls} with $\sourceG=0$.
  The main result of this section is stated as follows.

\begin{theorem}\label{thm.ls}
Fix $1\leq r\leq 2$, $m\ge 0$, and $\ell \in\R$.
Consider the Cauchy problem \eqref{ls} with $\sourceG=0$.  
For the soft potentials \eqref{hypSOFT} with $\wE\geq 0$, the solution of the
linearized homogeneous system satisfies the following time decay estimate
\begin{equation}
 \|w_\wN U(t)\solU_0\|_{L^2_p(\dot{H}^m_x)}
\lesssim
(1+t)^{-\sigma_{r,m}}\|w_{\wN+\wK}\solU_0\|_{L^2_pL^r_x} 
+ 
(1+t)^{-\wE/2}\|w_{\wN+\wE}\solU_0\|_{L^2_p(\dot{H}^m_x)},
 \label{thm.ls.1.soft}
\end{equation}
for any $\wK>2\sigma_{r,m}$.
Recall that  $\sigma_{r,m}$ is given by \eqref{rateLIN}.
\end{theorem}

We will prove this main theorem in the following few sub-sections.  In Section \ref{secL2} we develop the theory of the different components of the solution to the linear equation \eqref{rBoltz}.  The main point of this section is to derive a collection balance law equations and high-order moment equations for the coefficients of the different components of the linear solution to the relativistic Kinetic equation.
 Section \ref{sec.sub.tfli} then provides the relevant weighted instantaneous time-frequency Lyapunov inequality.
Finally  Section \ref{sec.tf} contains the proof of the time decay of solutions to the linear equation as stated in Theorem \ref{thm.ls} using the previous developments combined with an interpolation technique and the standard H{\"{o}}lder and Hausdorff-Young argument.

\subsection{Linear $L^2$ Bounds and Decay}\label{secL2}
We consider \eqref{ls} with $g=0$ as
\begin{gather}
 \partial_t f + \hat{p}\cdot \nabla_x f + L (f)
=
0,
\quad
f(0, x, p)=f_0(x,p),
\label{rBoltz}
\end{gather}
For the relativistic Maxwellian, $J$, we have the normalization $\int_{\R^3}J(p) dp=1$.  Initially, we introduce the notation for some integrals as follows
$$
\mu^0\eqdef \int_{\R^3}p^0 Jdp,  \quad 
\mu^{00}\eqdef \int_{\R^3}(p^0)^2 Jdp, \quad    
\mu^{11}\eqdef \int_{\R^3}p_1^2 Jdp, 
$$
$$\mu^{11}_0\eqdef \int_{\R^3}\frac{p_1^2}{p^0} Jdp,   \quad 
  \mu^{1122}_{00}\eqdef \int_{\R^3}\frac{p_1^2p_2^2}{(p^0)^2} Jdp,   \quad 
    \mu^{1111}_{00}\eqdef \int_{\R^3}\frac{p_1^4}{(p^0)^2} Jdp,
    $$
$$ 
\mu^{11}_{00}\eqdef \int_{\R^3}\frac{p_1^2}{(p^0)^2} Jdp.
$$
Note that the constants introduced above can be expressed in terms of the Bessel functions but, since we have no need to explore the delicate properties of Bessel functions, we will only utilize the expressions above.

We list some of the results in \cite{Strain2010}, which will be useful below.
\begin{lemma}
\label{nuEST}  
\cite{Strain2010}.
Consider  \eqref{nuDEF} with the soft potential kernel \eqref{hypSOFT}.  Then
$$
\nu(p) 
\approx
(p^0)^{-b/2}.
$$
More generally, 
$
\int_{\mathbb{R}^3} ~  dq 
\int_{\mathbb{S}^{2}} ~ d\omega
~ v_{\o} ~ \sigma(g,\theta)~ J^\alpha (q)
\approx
(p^0)^{-b/2}
$
for any $\alpha >0$.
\end{lemma}

Given a small $\epsilon>0$, choose a smooth cut-off function $\chi =\chi (g)$ satisfying
\begin{equation}
\chi (g)=
\left\{
\begin{array}{cl}
 1 & {\rm if } ~~  g \ge 2\epsilon
 \\
  0 &
  {\rm if } ~~ g \le \epsilon.
\end{array}
\right.
\label{cut}
\end{equation}
Now with \eqref{cut} and \eqref{compactK} we define
\begin{gather}
K_2^{1-\chi}(h) 
\eqdef
\int_{\mathbb{R}^3} ~  dq 
\int_{\mathbb{S}^{2}} ~ d\omega
~ \left( 1 - \chi(g) \right)
~ v_{\o} ~ \sigma(g,\theta)
\sqrt{J(q)} \sqrt{J(q^{\prime})} ~ h(p^{\prime })
\notag
\\
+
\int_{\mathbb{R}^3} ~  dq 
\int_{\mathbb{S}^{2}} ~ d\omega
~ \left( 1 - \chi(g) \right)
~ v_{\o} ~ \sigma(g,\theta)
\sqrt{J(q)}\sqrt{J(p^{\prime})} ~ h(q^{\prime}).
\label{kCUT}
\end{gather}
Define $K_1^{1-\chi}(h)$ similarly. We use the splitting
$
K
\eqdef
K^{1-\chi}
+
K^{\chi}$.  The following representation is derived in the Appendix of \cite{Strain2010}:
$$
K_i^\chi (h) = 
\int_{\mathbb{R}^3}dq~ k_i^\chi(p,q) ~ h(q),
\quad
i=1,2.
$$
Here is an estimate of $k_i^\chi(p,q)$:

\begin{lemma}
\label{boundK2}  
\cite{Strain2010}.
For the soft potentials \eqref{hypSOFT}, the kernel enjoys the estimate
$$
0\le k_2^{\chi}(p,q)
\le C_\chi \left( p^0 q^0 \right)^{-\zeta}
\left( p^0+ q^0 \right)^{-b/2}
e^{-c |p-q|},
\quad
C_\chi, c >0,
$$
with
$
\zeta
\eqdef
\min\left\{2-|\gamma|, 4-b,2\right\}/4
>0.
$
This estimate also holds for $k_1^\chi$.
\end{lemma}

We will also use the following estimate:
\begin{lemma}
\label{noKestimate}  
\cite{Strain2010}.
Fix any small $\eta >0,$ we may decompose  $K$ from \eqref{compactK} as 
$$
K = K_c + K_s.
$$
Here, for any $\ell\ge 0$, and for some $R=R(\eta)>0$ sufficiently large we have
$$
|\langle w^2_\ell  K_c h_1, h_2\rangle |
\le
C_{\eta} \|{\bf 1}_{\le R} h_1 \|_{L^2_p} \| {\bf 1}_{\le R} h_2 \|_{L^2_p},
$$
where ${\bf 1}_{\le R}$ is the indicator function of the ball of radius $R$ centered at zero.
We furthermore have the following estimate for the small part
\begin{equation}
|\langle w^2_{\ell}  K_s h_1, h_2\rangle |
\le 
 \eta 
\|  w_\ell h_1\nu^{\frac{1}{2}}\|_{L^2_p}
 \|  w_\ell h_2\nu^{\frac{1}{2}}\|_{L^2_p}.
\notag
\end{equation}
\end{lemma}
For every fixed $(t,x)$ the null
space of $L$ from \eqref{L} is given by the five dimensional space \cite{MR1379589}:
\begin{equation}\label{null}
{\mathcal N}\eqdef {\rm span}\left\{
\sqrt{J},
p_1 \sqrt{J}, p_2 \sqrt{J}, p_3\sqrt{J}, 
p^0\sqrt{J}\right\}.
\end{equation}
We define the orthogonal projection from $L^2(\mathbb{R}^3_p)$ onto the null space ${\mathcal N}$ by ${\bf P}$. 
Further expand ${\bf P} h$ as a linear combination of the basis in ~(\ref{null}): 
\begin{equation}
{\bf P} h
\eqdef
 \left\{a^h(t,x)+\sum_{j=1}^3 b_j^h(t,x)p_j+c^h(t,x)p^0\right\}\sqrt{J},
\label{hydro}
\end{equation}
where 
$$a^h=\int_{\R^3}h\sqrt{J}dp-\mu^0c^h,$$
$$b^h=\frac{\int_{\R^3}hp\sqrt{J}dp}{\mu^{11}},$$
$$c^h=\frac{\int_{\R^3}h(p^0\sqrt{J}-\mu^0\sqrt{J})dp}{\mu^{00}-(\mu^0)^2}.$$
We can then decompose  $f(t,x,p)$ as
\begin{equation}\label{E:f}
f={\bf P}f+\{{\bf I-P}\}f.
\end{equation}
With this decomposition we have the coercive estimate:

\begin{lemma}
\label{lowerN}  
\cite{Strain2010}.
$L\ge 0$.
$Lh = 0$ if and only if $h = {\bf P} h$.  
And
$\exists \delta_0>0$ such that
\begin{equation*}
\rmre \langle {L} h, h \rangle
\ge
\delta_0 
\|\nu^{1/2} \{ {\bf I - P } \} h \|_{L^2_p}^2.
\end{equation*}
\end{lemma}

We will see that for a solution to \eqref{rBoltz} the coefficients of ${\bf P}f$ satisfy balance laws:
\begin{gather}
\label{E:VPB3.9}
\partial_t a^f\left(1-\frac{(\mu^0)^2}{\mu^{00}}\right)
+\nabla_x\cdot b^f\left(\mu^{11}_0-\frac{\mu^{11}\mu^0}{\mu^{00}}\right)+\sum_{m=1}^3\partial_m \highB_m(\{{\bf I-P}\}f)=0,
\\
\label{E:VPB3.8}
\partial_t b_j^f\mu^{11}+\partial_j a^f\mu^{11}_0+\partial_j c^f\mu^{11}+\sum_{m=1}^3\partial_m \highG_{mj}(\{{\bf I-P}\}f)=0,
\quad 1 \le j \le 3, 
\\
\label{E:VPB3.9-2}
\partial_t c^f\left(\mu^0-\frac{\mu^{00}}{\mu^0}\right)
+\nabla_x\cdot b^f\left(\mu^{11}_0-\frac{\mu^{11}}{\mu^0}\right)+\sum_{m=1}^3\partial_m \highB_m(\{{\bf I-P}\}f)=0.
\end{gather}
These will be derived in the following developments.  

Now we are using the notation $\partial_m = \frac{\partial}{\partial x_m}$.
The high order moment functions, $\highG_{mj}(h)$, which are included in the above balance laws, are given by
\begin{equation}\label{E:A{mj}}
\highG_{mj}(h)\eqdef\int_{\R^3}\left(\frac{p_mp_j}{p^0}-\alpha_1\right)\sqrt{J(p)} h(p) dp,
\end{equation}
for $1\leq m,j\leq 3$, and $\alpha_1$ satisfies 
\begin{equation}\label{E:alpha1}
\frac{\mu^{11}}{\mu^0}(\mu^{11}_0-\alpha_1)-\mu^{1122}_{00}+\alpha_1\mu^{11}_0=0.
\end{equation}
Our choice of $\alpha_1$ is for a simpler expression in \eqref{E:VPB3.11}.  Furthermore $\highB_m(h)$ is
\begin{equation}\label{E:Bm}
\highB_m(h)\eqdef\int_{\R^3}p_m\left(\frac{1}{p^0}-\alpha_2\right)\sqrt{J(p)}h(p) dp
\end{equation}
for $1\leq m\leq 3$, and $\alpha_2=\frac{\mu^{11}_0}{\mu^{11}}.$ Our choice of $\alpha_2$ is for a simpler expression in \eqref{E:VPB3.13}.

We apply the decomposition (\ref{hydro}) and (\ref{E:f}) to the equation (\ref{rBoltz}). Multiplying (\ref{rBoltz}) by $\sqrt{J}$, $p_i\sqrt{J}$, $p^0\sqrt{J}$ for $1\leq i\leq 3$ and integrating over $\mathbb{R}^3$, one can get
\begin{eqnarray}
\partial_t\int_{\R^3}f\sqrt{J}dp+\int_{\R^3}\hat{p}\cdot\nabla_x f\sqrt{J}dp=0, \label{E:preVPB3.9}\\
\partial_t\int_{\R^3}f\sqrt{J}p_idp+\int_{\R^3}\hat{p}\cdot\nabla_x f\sqrt{J}p_idp=0, \quad 1\leq i\leq 3,
\label{E:preVPB3.8}
\\
\partial_t\int_{\R^3}f\sqrt{J}p^0dp+\int_{\R^3}\hat{p}\cdot\nabla_x f\sqrt{J}p^0dp=0, \label{E:preVPB3.7}
\end{eqnarray}
Using \eqref{E:f} and plugging (\ref{hydro}) into the above equation (\ref{E:preVPB3.7}) gives
\begin{equation}\label{E:VPB3.7}
\partial_t a^f\mu^0+\partial_t c^f\mu^{00}+\nabla_x\cdot b^f\mu^{11}=0.
\end{equation}
Plugging (\ref{hydro}) and \eqref{E:f} into the above equation (\ref{E:preVPB3.8}) gives \eqref{E:VPB3.8} using 
\eqref{E:VPB3.7}.

Plugging (\ref{hydro}) into the above equation (\ref{E:preVPB3.9}) gives
\begin{equation}\label{E:VPB3.9-0}
\partial_t a^f+\partial_t c^f\mu^0
+\nabla_x\cdot b^f\mu^{11}_0+\sum_{m=1}^3\partial_m \highB_m(\{{\bf I-P}\}f)=0,
\end{equation}
where recall $\highB_m(h)$ from \eqref{E:Bm}.
Then combining \eqref{E:VPB3.9-0} and (\ref{E:VPB3.7}) gives \eqref{E:VPB3.9}.
Similarly, combining \eqref{E:VPB3.9} with (\ref{E:VPB3.7}) grants \eqref{E:VPB3.9-2}.
This completes our derivation of the balance laws.  

To continue, we rewrite \eqref{rBoltz} as
\begin{equation}\label{E:rBoltz2}
\partial_t f+\hat{p}\cdot \nabla_x ({\bf P}f)=-L(\{{\bf I-P}\}f)-\hat{p}\cdot \nabla_x(\{{\bf I-P}\}f)\eqdef R.
\end{equation}
Using \eqref{E:rBoltz2}, we get
\begin{align}
\partial_t \highG_{jj}(\{{\bf I-P}\}f)&=\highG_{jj}(R-\hat{p}\cdot \nabla_x ({\bf P}f)-\partial_t ({\bf P}f))\notag\\
&=\highG_{jj}(R)-\int_{\R^3}\left(\frac{p_j^2}{p^0}-\alpha_1\right)\sqrt{J}\left(\hat{p}\cdot \nabla_x ({\bf P}f)+\partial_t ({\bf P}f)\right)\,dp\notag\\
&=\highG_{jj}(R)-\sum_{k\neq j}\partial_k b_k^f(\mu^{1122}_{00}-\alpha_1\mu^{11}_0)-\partial_j b_j^f(\mu^{1111}_{00}-\alpha_1\mu^{11}_0)\label{E:VPB3.11-0}\\
&\quad -\partial_t a^f(\mu^{11}_0-\alpha_1)-\partial_t c^f(\mu^{11}-\alpha_1\mu^0).\notag
\end{align}
From \eqref{E:VPB3.7}, we have
\begin{equation}\label{E:VPB3.11-1}
\partial_t a^f(\mu^{11}_0-\alpha_1)+\partial_t c^f\frac{\mu^{00}(\mu^{11}_0-\alpha_1)}{\mu^0}+\nabla_x\cdot b^f\frac{\mu^{11}(\mu^{11}_0-\alpha_1)}{\mu^0}=0.
\end{equation}
Combining \eqref{E:VPB3.11-0}, \eqref{E:VPB3.11-1} and \eqref{E:alpha1}, we get
\begin{multline}
\partial_t \highG_{jj}(\{{\bf I-P}\}f)=\highG_{jj}(R)+\partial_t c^f\left(\frac{\mu^{00}}{\mu^0}(\mu^{11}_0-\alpha_1)-\mu^{11}+\alpha_1\mu^0\right)\\
+\partial_j b_j^f\left(\mu^{1122}_{00}-\mu^{1111}_{00}\right).
\label{E:VPB3.11}
\end{multline}
For $j\neq m$, using \eqref{E:rBoltz2}, we get
\begin{align}
\partial_t \highG_{mj}(\{{\bf I-P}\}f)&=\highG_{mj}(R-\hat{p}\cdot \nabla_x ({\bf P}f)-\partial_t ({\bf P}f))\notag\\
&=\highG_{mj}(R)-\int_{\R^3}\left(\frac{p_jp_m}{p^0}-\alpha_1\right)\sqrt{J}\left(\hat{p}\cdot \nabla_x ({\bf P}f)+\partial_t ({\bf P}f)\right)\,dp\notag\\
&=\highG_{mj}(R)-(\partial_m b_j^f+\partial_j b^f_m)\mu^{1122}_{00}-\alpha_1\sum_{k=1}^3 \partial_k \Lambda_k (\{{\bf I-P}\}f).\label{E:VPB3.12}
\end{align}
Again using \eqref{E:rBoltz2}, \eqref{E:preVPB3.7} and our choice of $\alpha_2$, we get
\begin{align}
\partial_t \highB_{m}(\{{\bf I-P}\}f)&=\highB_{m}(R-\hat{p}\cdot \nabla_x ({\bf P}f)-\partial_t ({\bf P}f))\notag\\
&=\highB_{m}(R)-\int_{\R^3}p_m\left(\frac{1}{p^0}-\alpha_2\right)\sqrt{J}\left(\hat{p}\cdot \nabla_x ({\bf P}f)+\partial_t ({\bf P}f)\right)\,dp\notag\\
&=\highB_{m}(R)-\partial_m a^f(\mu^{11}_{00}-\alpha_2\mu^{11}_0)-\partial_m c^f(\mu^{11}_{0}-\alpha_2\mu^{11})\notag\\
&\quad-\partial_t b_m^f(\mu^{11}_{0}-\alpha_2\mu^{11})\notag\\
&=\highB_m(R)-\partial_m a^f(\mu^{11}_{00}-\alpha_2\mu^{11}_0).\label{E:VPB3.13}
\end{align}
By \eqref{E:A{mj}} and \eqref{null}, we have
\begin{align}\label{E:VPB3.15-a}
\sum_{j=1}^3 \partial_t \partial_m \highG_{jj}(\{{\bf I-P}\}f)&=\int_{\R^3}\left(\frac{|p|^2}{p^0}-3\alpha_1\right)\sqrt{J}\partial_t \partial_m(\{{\bf I-P}\}f)dp\notag\\
&=\int_{\R^3}\frac{-1}{p^0}\sqrt{J}(\{{\bf I-P}\}\partial_t \partial_m f)dp.
\end{align}
With this calculation in mind, we introduce the definition:
$$
A(h)\eqdef\int_{\R^3}\frac{\sqrt{J(p)}}{p^0}h(p)dp.
$$
On the other hand, by \eqref{E:VPB3.11}, we have
\begin{align}
&\quad\sum_{j=1}^3 \partial_t \partial_m \highG_{jj}(\{{\bf I-P}\}f)\notag\\
&=\sum_{j=1}^3 \partial_m \highG_{jj}(R)+3\partial_m \partial_t c^f\left(\frac{\mu^{00}}{\mu^0}(\mu^{11}_0-\alpha_1)-\mu^{11}+\alpha_1\mu^0\right)\notag\\
&\quad+\partial_m\nabla_x\cdot b^f\left(\mu^{1122}_{00}-\mu^{1111}_{00}\right).\label{E:VPB3.15-a'}
\end{align}
Combining \eqref{E:VPB3.15-a} and \eqref{E:VPB3.15-a'} gives
\begin{multline}\label{E:VPB3.15-b}
\sum_{j=1}^3 \partial_m \highG_{jj}(R)+3\partial_m \partial_t c^f\left(\frac{\mu^{00}}{\mu^0}(\mu^{11}_0-\alpha_1)-\mu^{11}+\alpha_1\mu^0\right)\\
+\partial_m\nabla_x\cdot b^f\left(\mu^{1122}_{00}-\mu^{1111}_{00}\right)
= 
- A\left(\{{\bf I-P}\}\partial_t \partial_m f\right).
\end{multline}
By \eqref{E:VPB3.11} and \eqref{E:VPB3.12}, for any constant $\beta$ we have
\begin{eqnarray}\label{E:VPB3.15-c}
&&\partial_t\sum_{j=1}^3\partial_j \highG_{jm}(\{{\bf I-P}\}f)+\beta\partial_t\partial_m \highG_{mm}(\{{\bf I-P}\}f)\notag\\
&&=\sum_{j=1}^3\partial_j \highG_{jm}(R)+\beta\partial_m \highG_{mm}(R)-\Delta_x b^f_m\mu^{1122}_{00}-\partial_m\nabla_x\cdot b^f\mu^{1122}_{00}\notag\\
&& +(\beta+1)\partial_m\partial_t c^f\left(\frac{\mu^{00}}{\mu^0}(\mu^{11}_0-\alpha_1)-\mu^{11}+\alpha_1\mu^0\right)\notag\\
&& +(\beta+1)\partial_m\partial_m b^f_m\left(\mu^{1122}_{00}-\mu^{1111}_{00}\right)\notag\\
&&+2\partial_m\partial_m b^f_m\mu^{1122}_{00}-\alpha_1\sum_{j\neq m}\sum_{k=1}^3\partial_j \partial_k \Lambda_k (\{{\bf I-P}\}f).
\end{eqnarray}
Choose $\beta$ such that $$\frac{\beta+1}{3}=\frac{-\mu^{1122}_{00}}{\mu^{1122}_{00}-\mu^{1111}_{00}}.$$ Then \eqref{E:VPB3.15-c} becomes
\begin{eqnarray}\label{E:VPB3.15-c'}
&&\partial_t\sum_{j=1}^3\partial_j \highG_{jm}(\{{\bf I-P}\}f)+\beta\partial_t\partial_m \highG_{mm}(\{{\bf I-P}\}f)\notag\\
&&=\sum_{j=1}^3\partial_j \highG_{jm}(R)+\beta\partial_m \highG_{mm}(R)-\Delta_x b^f_m\mu^{1122}_{00}-\partial_m\nabla_x\cdot b^f\mu^{1122}_{00}\notag\\
&& +(\beta+1)\partial_m\partial_t c^f\left(\frac{\mu^{00}}{\mu^0}(\mu^{11}_0-\alpha_1)-\mu^{11}+\alpha_1\mu^0\right)\notag\\
&&-\partial_m\partial_m b^f_m\mu^{1122}_{00}-\alpha_1\sum_{j\neq m}\sum_{k=1}^3\partial_j \partial_k \Lambda_k (\{{\bf I-P}\}f).
\end{eqnarray}
Also \eqref{E:VPB3.15-b} implies
\begin{multline}\label{E:VPB3.15-b'}
\frac{\beta+1}{3}\sum_{j=1}^3 \partial_m \highG_{jj}(R)+(\beta+1)\partial_m \partial_t c^f\left(\frac{\mu^{00}}{\mu^0}(\mu^{11}_0-\alpha_1)-\mu^{11}+\alpha_1\mu^0\right)\\
-\partial_m\nabla_x\cdot b^f\mu^{1122}_{00}
= -\frac{\beta+1}{3}
A\left(\{{\bf I-P}\}\partial_t \partial_m f\right).
\end{multline}
Combining \eqref{E:VPB3.15-b'} and \eqref{E:VPB3.15-c'} gives
\begin{eqnarray}\label{E:VPB3.15}
&&-\partial_t\left(\sum_{j =1}^3\partial_j \highG_{jm}(\{{\bf I-P}\}f)+\beta\partial_m \highG_{mm}(\{{\bf I-P}\}f)\right)\notag
\\
&&
- \frac{\beta+1}{3}
\partial_t A\left(\{{\bf I-P}\}\partial_m f\right)
\notag\\
&&-\Delta_x b^f_m\mu^{1122}_{00}-\partial_m\partial_m b^f_m\mu^{1122}_{00}\notag\\
&&=\frac{\beta+1}{3}\sum_{j=1}^3\partial_m \highG_{jj}(R)-\sum_{j=1}^3\partial_j \highG_{mj}(R)-\beta\partial_m \highG_{mm}(R)\\
&&+\alpha_1\sum_{j\neq m}\sum_{k=1}^3\partial_j \partial_k \Lambda_k (\{{\bf I-P}\}f).
\notag
\end{eqnarray}
We are now ready to prove the following lemma.  This lemma is a key step in creating a time-frequency Lyapunov function later on.

\begin{lemma}\label{lem.li.fr}
There is a free energy functional $\CE_{free}(\widehat{f}(t,k))$ which
is local in the time and frequency variables, and takes the form of
\begin{eqnarray}
 \CE_{free}(\widehat{f}(t,k))  
&\eqdef & \kappa_1 \sum_m \left(\sum_j
  \frac{i k_j}{1+|k|^2}
\highG_{jm}(\{\FI-\FP\}\widehat{f}) \mid \, -b^{\widehat{f}}_m\right)
\nonumber\\
&&+
\kappa_1 \sum_m \left(\beta \frac{i k_m}{1+|k|^2}
\highG_{mm}(\{\FI-\FP\}\widehat{f}) \mid \, -b^{\widehat{f}}_m\right)
\nonumber\\
&&\ \ \ 
+
\kappa_1 \sum_m \left(\frac{\beta+1}{3}\frac{i k_m}{1+|k|^2}A(\{\FI-\FP\}\widehat{f})\mid \, -b^{\widehat{f}}_m\right)\nonumber\\
&&\ \ \ +\kappa_1 \sum_{j}\left( \highB_j (\{\FI-\FP\}\widehat{f}) \mid
\,\frac{i k_j}{1+|k|^2} a^{\widehat{f}}\right)
\nonumber\\
&&\ \ \ 
+\sum_m \left( b^{\widehat{f}}_m \mid\, \frac{ i k_m}{1+|k|^2}
(\mu^{11}_0a^{\widehat{f}}+\mu^{11}c^{\widehat{f}})\right),
\label{E:VPB3.16}
\end{eqnarray}
for some constant $\kappa_1>0$, such that one has
\begin{eqnarray}
&\dis \pa_t \rmre\,  \CE_{free}(\widehat{f}(t,k))+ \la
\frac{|k|^2}{1+|k|^2}\left(\left| a^{\widehat{f}} \right|^2+\left| b^{\widehat{f}} \right|^2+\left| c^{\widehat{f}} \right|^2
\right)\nonumber\\
&\dis \leq C\|\{\FI-\FP\}\widehat{f}\nu^{\frac{1}{2}}\|_{L^2_p}^2
\label{E:VPB3.17}
\end{eqnarray}
for any $t\geq 0$ and $k\in\R^3$.
\end{lemma}

\begin{proof}
The proof of this Lemma 3.1 uses the overall strategy in \cite{MR2754344}, even though significant new difficulties arise because we are incorporating the effects of special relativity.  
We shall make estimates on $b^{\widehat{f}}$, $a^{\widehat{f}}$ and
$\mu^{11}_0a^{\widehat{f}}+\mu^{11}c^{\widehat{f}}$ individually and then take the
proper linear combination to deduce the desired free energy
inequality \eqref{E:VPB3.17}. Firstly, notice that
\begin{equation*}
  \CF a^f=a^{\CF f} =a^{\widehat{f}},\ \  \CF b^f=b^{\CF f}=b^{\widehat{f}}, \ \ \CF c^f=c^{\CF f} =c^{\widehat{f}},
\end{equation*}
and likewise for the high-order moment functions $\highG_{jm}(\cdot)$, $\highB_j(\cdot)$ and $A(\cdot)$.

\medskip
\noindent{\it Estimate on $b^{\widehat{f}}$.}  We claim that for any
$0<\de_1<\mu^{1122}_{00}$, it holds that
\begin{eqnarray} \notag
&&\pa_t {\rmre}\sum_m \left(\sum_j i k_j
\highG_{jm}(\{\FI-\FP\}\widehat{f})+\beta i k_m
\highG_{mm}(\{\FI-\FP\}\widehat{f})\right.\\
&&\left.+\frac{\beta+1}{3}i k_mA(\{\FI-\FP\}\widehat{f})\mid \,-b^{\widehat{f}}_m\right)+(\mu^{1122}_{00}-\de_1)|k|^2\left| b^{\widehat{f}} \right|^2
\nonumber \\
&&\leq \de_1
|k|^2
\left|\mu^{11}_0a^{\widehat{f}}+\mu^{11}c^{\widehat{f}}\right|^2
+
\de_1|k|^2\left| a^{\widehat{f}} \right|^2\nonumber \\
&&\ \ \
+\frac{C}{\de_1}(1+|k|^2)\|\{\FI-\FP\}\widehat{f}\nu^{\frac{1}{2}}\|_{L^2_p}^2.
\label{E:VPB3.18}
\end{eqnarray}
In fact, the Fourier transform of \eqref{E:VPB3.15} gives
\begin{eqnarray*}
 &&\dis-\pa_t\left[\sum_{j}i k_j \highG_{jm}(\{\FI-\FP\}\widehat{f})
 +\beta i
 k_m \highG_{mm}(\{\FI-\FP\}\widehat{f})\right.\\
 &&\left.+\frac{\beta+1}{3}i k_mA(\{\FI-\FP\}\widehat{f})\right]+(|k|^2b_m^{\widehat{f}}
 +k_m^2
 b_m^{\widehat{f}})\mu^{1122}_{00}\nonumber\\
 &&\dis =\alpha_1\sum_{j\neq m}\sum_{l=1}^3k_j k_l \Lambda_l (\{{\bf I-P}\}f)\\ 
 &&+\frac{\beta+1}{3}\sum_{j\neq m}i k_m \highG_{jj}(\widehat{R})-\sum_{j}i
 k_j
 \highG_{jm}(\widehat{R})-\beta i k_m\highG_{mm}(\widehat{R}).
\end{eqnarray*}
Taking further the complex inner product with $b_m^{\widehat{f}}$
gives
\begin{eqnarray}
 &&\dis\pa_t\left(\sum_{j}i k_j \highG_{jm}(\{\FI-\FP\}\widehat{f})
 +\beta i
 k_m \highG_{mm}(\{\FI-\FP\}\widehat{f})\right.\nonumber\\
 &&
 \left.+\frac{\beta+1}{3}i k_mA(\{\FI-\FP\}\widehat{f}) \mid \,-b_m^{\widehat{u}}\right)
 +(|k|^2+k_m^2)
 \left|b_m^{\widehat{f}}\right|^2\mu^{1122}_{00}
 \nonumber\\
 &&\dis =I_1+I_2,
 \label{E:VPB3.19}
\end{eqnarray}
where
\begin{eqnarray}
I_1 & = & \left( \alpha_1\sum_{j\neq m}\sum_{l=1}^3k_j k_l \Lambda_l (\{{\bf I-P}\}f)+\frac{\beta+1}{3}\sum_{j\neq m}i k_m \highG_{jj}(\widehat{R})\right.
 \nonumber\\
 &&\ \ \left.-\sum_{j}i
 k_j
 \highG_{jm}(\widehat{R})-\beta i k_m\highG_{mm}(\widehat{R})\mid \, b_m^{\widehat{u}}\right),  
 \nonumber\\
 I_2 &=&
 \left(\sum_{j}i k_j \highG_{jm}(\{\FI-\FP\}\widehat{f})
 +\beta i
 k_m \highG_{mm}(\{\FI-\FP\}\widehat{f})\right.
 \nonumber\\
 &&\ \ \left.+\frac{\beta+1}{3}i k_mA(\{\FI-\FP\}\widehat{f})\mid  \,
 -\pa_tb_m^{\widehat{f}}\right).
 \notag
\end{eqnarray}
Then $I_1$ is bounded by
\begin{eqnarray*}
  \left| I_1 \right| &\leq & \de_1 |k|^2
  \left|b_m^{\widehat{f}}\right|^2
  +
  \frac{C}{\de_1}\sum_{jm} \left|\highG_{jm}(\widehat{R})\right|^2
  +
  \frac{C}{\delta_1}
  \sum_{l=1}^3(1+|k|^2) \left|\Lambda_l (\{{\bf I-P}\}f)\right|^2.
\end{eqnarray*}
For $I_2$, one can use the Fourier transforms of \eqref{E:VPB3.8}:
\begin{equation}
 \pa_t b^{\widehat{f}}_j \mu^{11} +i k_j
(a^{\widehat{f}}\mu^{11}_0+c^{\widehat{f}}\mu^{11})
+\sum_{m}i k_m \highG_{jm}(\{\FI- \FP\}\widehat{f})=0,
\label{E:VPB3.20}
\end{equation}
to estimate it as 
\begin{multline*}
  \left| I_2 \right| \leq \de_1 |k|^2
  \left|\mu^{11}_0a^{\widehat{f}}+\mu^{11}c^{\widehat{f}}\right|^2
  \\
  +\frac{C}{\de_1}|k|^2
  \left(
   \sum_{j,m} \left|\highG_{jm}(\{\FI- \FP\}\widehat{f})\right|^2+\left|A(\{\FI- \FP\}\widehat{f})\right|^2
  \right).
\end{multline*}
On the other hand, notice from \eqref{E:rBoltz2} that
\begin{equation*}
\widehat{R}=-\hat{p}\cdot ik\{\FI-\FP\}\widehat{f}+L
\{\FI-\FP\}\widehat{f},
\end{equation*}
which implies
\begin{eqnarray*}
 \sum_{j,m} \left|\highG_{jm}(\widehat{R})\right|^2+\left|A(\widehat{R})\right|^2
 &\leq &
 C(1+|k|^2)\|\{\FI-\FP\}\widehat{f}\nu^{\frac{1}{2}}\|_{L^2_p}^2.
\end{eqnarray*}
Thus we see that \eqref{E:VPB3.18} follows from taking the real part of
\eqref{E:VPB3.19} and plugging the estimates for $I_1$ and $I_2$ into
that.
\medskip

\noindent{\it Estimate on $a^{\widehat{f}}$.}  We claim that for any
$0<\de_2<\mu^{11}_{00}-\alpha_2\mu^{11}_0$, it holds that
\begin{eqnarray}
&&\pa_t \rmre \sum_{j}\left( \highB_j (\{\FI-\FP\}\widehat{f}) \mid \, i k_j
a^{\widehat{f}}\right)+(\mu^{11}_{00}-\alpha_2\mu^{11}_0-\de_2)|k|^2 \left| a^{\widehat{f}} \right|^2  \nonumber\\
&& \leq \de_2
|k|^2\left| b^{\widehat{f}} \right|^2+\frac{C}{\de_2}(1+|k|^2)\|\{\FI-\FP\}\widehat{f}\nu^{\frac{1}{2}}\|_{L^2_p}^2.\label{E:VPB3.22}
\end{eqnarray}
In fact, similarly as before, from the Fourier transform of
\eqref{E:VPB3.13}
\begin{equation*}
\pa_t \highB_j(\{\FI- \FP\}\widehat{f})+i k_j
a^{\widehat{f}}(\mu^{11}_{00}-\alpha_2\mu^{11}_0)=\highB_j(\widehat{R}),
\end{equation*}
one can get
\begin{eqnarray}
&\dis \pa_t \left(\highB_j(\{\FI- \FP\}\widehat{f}) \mid \, i k_j
a^{\widehat{f}}\right) +|k_j|^2
\left| a^{\widehat{f}} \right|^2(\mu^{11}_{00}-\alpha_2\mu^{11}_0)  \nonumber\\
&\dis =
\left(\highB_j(\widehat{R}) \mid \, i k_j
a^{\widehat{f}}\right)
+
\left( \highB_j(\{\FI- \FP\}\widehat{f}) \mid \, i k_j
\pa_t a^{\widehat{f}}\right)
=I_3+I_4.\label{E:VPB3.23}
\end{eqnarray}
Now $I_3$ is bounded by  
\begin{eqnarray*}
  \left| I_3 \right|  &\leq & \de_2|k_j|^2 \left| a^{\widehat{f}} \right|^2
  +
  \frac{C}{\de_2}\sum_j 
  \left|\highB_j(\widehat{R})\right|^2,
\end{eqnarray*}
and from the Fourier transform of \eqref{E:VPB3.9}
\begin{eqnarray*}
    &&\pa_t a^{\widehat{f}}+i k \cdot b^{\widehat{f}}\frac{\mu^{11}_{0}\mu^{00}-\mu^{11}\mu^0}{\mu^{00}-(\mu^0)^2}+\sum_{j}i
    k_j
\highB_j(\{\FI- \FP\}\widehat{f})\frac{\mu^{00}}{\mu^{00}-(\mu^0)^2}=0,
\end{eqnarray*}
where $I_4$ is bounded by
\begin{equation*}
    \left| I_4 \right| \leq\de_2|k|^2  
    \left| b^{\widehat{f}} \right|^2
    +\frac{C}{\de_2}\left( 1+|k|^2 \right) \sum_j \left|\highB_j(\{\FI- \FP\}\widehat{f})\right|^2.
\end{equation*}
Notice that similar to $\highG_{jm}$, it holds that
\begin{equation*}
 \left|\highB_{j}(\widehat{R})\right|^2
 \leq
 C(1+|k|^2)\|\{\FI-\FP\}\widehat{f}\nu^{\frac{1}{2}}\|_{L^2_p}^2.
\end{equation*}
Then, \eqref{E:VPB3.22} follows from \eqref{E:VPB3.23} by
taking summation over  $i$, taking the real part and then
applying the estimates of $I_3$ and $I_4$.
\medskip

\noindent{\it Estimate on $\mu^{11}_0a^{\widehat{f}}+\mu^{11}c^{\widehat{f}}$.} We notice that $c^{\widehat{f}}$ is a linear combination of $\mu^{11}_0a^{\widehat{f}}+\mu^{11}c^{\widehat{f}}$ and $a^{\widehat{f}}$. So our estimates on $\mu^{11}_0a^{\widehat{f}}+\mu^{11}c^{\widehat{f}}$ and $a^{\widehat{f}}$ imply the estimate of $c^{\widehat{f}}$. We
claim that for any $0<\de_3<\frac{1}{\mu^{11}}$, it holds that
\begin{multline}
\pa_t \rmre 
\sum_m \left( b^{\widehat{f}}_m \mid \, i k_m (\mu^{11}_0a^{\widehat{f}}+\mu^{11}c^{\widehat{f}})\right)
+\left(\frac{1}{\mu^{11}}-\de_3\right)|k|^2
\left|\mu^{11}_0a^{\widehat{f}}+\mu^{11}c^{\widehat{f}}\right|^2
\\
\leq
C|k|^2\left| b^{\widehat{f}} \right|^2+\frac{C}{\de_3}|k|^2\|\{\FI-\FP\}\widehat{f}\nu^{\frac{1}{2}}\|_{L^2_p}^2.
\label{E:VPB3.24}
\end{multline}
In fact, by taking the complex inner product with $i k_j
(\mu^{11}_0a^{\widehat{f}}+\mu^{11}c^{\widehat{f}}) $ and then taking summation over
$1\leq j\leq n$, it follows from \eqref{E:VPB3.20} that
\begin{eqnarray}
 &&\dis \pa_t \sum_j\left( b^{\widehat{f}}_j \mid \, i k_j
(\mu^{11}_0a^{\widehat{f}}+\mu^{11}c^{\widehat{f}}) \right) 
+\frac{|k|^2}{\mu^{11}}
\left|\mu^{11}_0a^{\widehat{f}}+\mu^{11}c^{\widehat{f}}\right|^2 \nonumber\\
&&\dis= 
\sum_{j,m}
\left( \frac{-i
k_m}{\mu^{11}}\highG_{jm}(\{\FI-\FP\}\widehat{f}) \mid \, i k_j
(\mu^{11}_0a^{\widehat{f}}+\mu^{11}c^{\widehat{f}})\right)  
\nonumber\\
&&\ \ \ + \sum_j\left( b^{\widehat{f}} \mid \, i k_j \pa_t
(\mu^{11}_0a^{\widehat{f}}+\mu^{11}c^{\widehat{f}}) \right)
\nonumber\\
&&=I_5+I_6.\label{E:VPB3.25}
\end{eqnarray}
Now $I_5$ is bounded similar to previous estimates as
\begin{eqnarray*}
\left| I_5 \right| 
\leq 
\de_3|k|^2 \left| \mu^{11}_0a^{\widehat{f}}+\mu^{11}c^{\widehat{f}} \right|^2+\frac{C}{\de_3}|k|^2\|\{\FI-\FP\}\widehat{f}\nu^{\frac{1}{2}}\|_{L^2_p}^2.
\end{eqnarray*}
For $I_6$, it holds that
\begin{eqnarray*}
\left| I_6 \right| 
\leq 
C|k|^2 \left| b^{\widehat{f}} \right|^2+C|k|^2\|\{\FI-\FP\}\widehat{f}\nu^{\frac{1}{2}}\|_{L^2_p}^2,
\end{eqnarray*}
where we have used the Fourier transform of \eqref{E:VPB3.9} as
\begin{eqnarray*}
    &&\pa_t a^{\widehat{f}}+i k \cdot b^{\widehat{f}}\frac{\mu^{11}_{0}\mu^{00}-\mu^{11}\mu^0}{\mu^{00}-(\mu^0)^2}+\sum_{j}i
    k_j
\highB_j(\{\FI- \FP\}\widehat{f})\frac{\mu^{00}}{\mu^{00}-(\mu^0)^2}=0,
\end{eqnarray*}
and the Fourier transform of \eqref{E:VPB3.9-2}
\begin{eqnarray*}
    &&\pa_t c^{\widehat{f}}+i k \cdot b^{\widehat{f}}\frac{\mu^{11}_{0}\mu^{0}-\mu^{11}}{(\mu^0)^2-\mu^{00}}+\sum_{j}i
    k_j
\highB_j(\{\FI- \FP\}\widehat{f})\frac{\mu^{0}}{(\mu^0)^2-\mu^{00}}=0.
\end{eqnarray*}
Therefore, \eqref{E:VPB3.17} follows from the proper
linear combination of \eqref{E:VPB3.18}, \eqref{E:VPB3.22}
and \eqref{E:VPB3.24} by taking $0<\de_1,\de_2,\de_3<1$ small
enough and also $\kappa_0>0$ large enough. This completes the proof
of Lemma \ref{lem.li.fr}.
\end{proof}

\subsection{Weighted time-frequency Lyapunov inequality}\label{sec.sub.tfli}
In this subsection, we shall construct the desired time-frequency
Lyapunov functional.

\subsubsection{Estimate on the microscopic dissipation}

The first step in our construction of the time-frequency Lyapunov
functional is to estimate the microscopic dissipation on the basis of the
coercivity property in Lemma \ref{lowerN} of $\FL$. 

Consider  \eqref{rBoltz}, taking the Fourier transform in $x$ grants us
\begin{equation}\label{ls-1f}
  \dis     \pa_t \hat{\solU}+\rmi\vel\cdot k ~ \hat{\solU}+\FL \hat{\solU} =0.
\end{equation}
Then we multiply equation \eqref{ls-1f} with $\overline{\hat{\solU}}(t,k)$ and integrate over $\R^3_\vel$ to achieve
\begin{equation}
\notag
    \frac{1}{2}\pa_t \|\hat{\solU}\|_{L^2_p}^2+\rmre \ang{\FL \hat{\solU}, \hat{\solU} }
    =0.
\end{equation}
By Lemma \ref{lowerN}, one has that
\begin{equation}\label{diss-micr}
\pa_t \|\hat{\solU}\|_{L^2_p}^2
+
\la \| \nu^{1/2} \{\FI-\FP\} \hat{\solU}\|_{L^2_p}^2
\lesssim
0.
\end{equation}
This is the first main estimate which we will use in the following.

\subsubsection{Macroscopic time-frequency weighted inequality}
In this section we prove the following instantaneous Lyapunov inequality
with a velocity weight $\wN \in \R$:
\begin{multline}\label{macroWeightINEQ}
\frac{1}{2}\frac{d}{dt}\| w_{\wN}\{\FI-\FP\}\hat{\solU}\|^2_{L^2_p} +
\la \| \nu^{1/2}w_{\wN}\{\FI-\FP\}\hat{\solU}\|^2_{L^2_{p}}
\\
\le C_\lambda |k|^2 \|\nu^{1/2} \hat{\solU}\|^2_{L^2_{p}}
+
 C \|{\bf 1}_{\leq C}\{\FI-\FP\}\hat{\solU}\|_{L^2_p}^2.
\end{multline}
We split the solution $\solU$ to equation \eqref{rBoltz} into $\solU=\FP \solU + \{\FI-\FP\}\solU$, take the Fourier transform as in \eqref{ls-1f},  and then apply $\{\FI-\FP\}$ to the resulting equation:
\begin{multline}\notag
\pa_t \{\FI-\FP\}\hat{\solU} + \rmi \vel \cdot k \{\FI-\FP\}\hat{\solU} 
+
\FL\{\FI-\FP\} \hat{\solU}\\
=-\{\FI-\FP\}(\rmi \vel \cdot k \FP \hat{\solU} )
 +\FP (\rmi \vel \cdot k \{\FI-\FP\}\hat{\solU}).
\end{multline}
Multiply the last equation by $w_{2\wN}\{\FI-\FP\}\overline{\hat{\solU}}$ and integrate in $\R^3_\vel$ to obtain
\begin{equation}
\label{app.vw.p05.int1.new}
\frac{1}{2}\frac{d}{dt}\| w_{\wN}\{\FI-\FP\}\hat{\solU}\|^2_{L^2_p} +
 \rmre \langle w_{2\wN}\FL\{\FI-\FP\} \hat{\solU},\{\FI-\FP\} \hat{\solU}\rangle
 = \Ga_1,
\end{equation}
where 
\begin{multline*}
 \Ga_1 =
-\rmre\left\langle
\{\FI-\FP\}(\rmi \vel \cdot k\FP \hat{\solU} ),
w_{2\wN}  \{\FI-\FP\}\hat{\solU}
\right\rangle
\\
 +\rmre \left\langle
 \FP (\rmi \vel \cdot k \{\FI-\FP\}\hat{\solU}),
w_{2\wN}  \{\FI-\FP\}\hat{\solU}  \right\rangle.
\end{multline*}
As a result of the rapid decay in the coefficients of \eqref{hydro} we obtain
$$
\left| \Ga_1 \right| \le
\eta \|\nu^{1/2} w_{\wN}\{\FI-\FP\}\hat{\solU}\|^2_{L^2_{p}}
+
C_\eta |k|^2 \left( \| w_{-\wE}\{\FI-\FP\}\hat{\solU}\|^2_{L^2_p}+\|\nu^{1/2} \FP \hat{\solU}\|^2_{L^2_p} \right), 
$$
which holds for any small $\eta >0$ and any large $\wE>0$. 

By \eqref{L} and Lemma \ref{noKestimate}, we have 
$$
\rmre  \langle w_{2\wN}\FL\{\FI-\FP\} \hat{\solU},\{\FI-\FP\} \hat{\solU}\rangle
 \ge \la \| w_{\wN}\{\FI-\FP\}\hat{\solU}\nu^{\frac{1}{2}}\|_{L^2_p}^2
 -
 C  \|{\bf 1}_{\leq C}\{\FI-\FP\}\hat{\solU}\|_{L^2_p}^2.
$$
This holds for a small $\la >0$ and a large $C>0$.  (We have used the fact that 
Lemma \ref{noKestimate} indeed holds for any $\ell \in \R$, as can be seen from the proofs in
\cite[Lemma 3.3]{Strain2010} and \cite[Lemma 3.6]{Strain2010}.)  Plugging the last few estimates into \eqref{app.vw.p05.int1.new} proves \eqref{macroWeightINEQ}.

We furthermore remark, following the same procedure as above, that we obtain
\begin{equation}\label{macroWeightINEQ.noM}
\frac{1}{2}\frac{d}{dt}\| w_{\wN}\hat{\solU}\|^2_{L^2_p} +
\la \| \nu^{1/2}w_{\wN}\hat{\solU}\|^2_{L^2_{p}}
\le 
 C  \|{\bf 1}_{\leq C}\hat{\solU}\|_{L^2_p}^2.
\end{equation}
In other words, if we multiply \eqref{ls-1f} by $w^{2\wN}\overline{\hat{\solU}}(t,k)$, integrate in $\R^3_\vel$ and use the same estimates as in the last case it follows that we obtain \eqref{macroWeightINEQ.noM}.

\subsubsection{Derivation of  time-frequency Lyapunov inequality}
Now we prove  

\begin{theorem}\label{thm.tfli}
Fix $\wN\in \R$.
Let $\solU$ be the solution to the Cauchy problem \eqref{ls}
with $g=0$. Then there is a time-frequency functional $\CE_\wN(t,k)$
such that
\begin{equation}\label{thm.tfli.1}
    \CE_\wN(t,k) \approx \| w_\ell \hat{\solU}\|_{L^2_p}^2,
\end{equation}
where for any $t\geq 0$ and $k\in \threed$ we have 
\begin{equation}\label{thm.tfli.2}
\dis\pa_t \CE_{\wN} (t,k)+\la \left( 1 \wedge |k|^2\right) \| \nu^{1/2}w_{\wN}  \hat{\solU}\|_{L^2_{p}}^2
\lesssim
 0.
\end{equation}
Here $1 \wedge |k|^2 \eqdef\min\{1,|k|^2\}$.
\end{theorem}

\begin{proof}
We first define 
\begin{equation}\label{thm.tfli.p1}
    \CE(t,k)\eqdef \|\hat{\solU}\|_{L^2_p}^2
    +\kappa_3 \CE_{free}(t,k),
\end{equation}
for a constant $\kappa_3>0$ to be determined later, where $\CE_{free}(t,k)$ is given in \eqref{E:VPB3.16}. 
One can fix $\kappa_3>0$ small enough such that 
$
\CE(t,k)\approx \|\hat{\solU}\|_{L^2_p}^2.
$

A linear combination of \eqref{diss-micr} and \eqref{E:VPB3.17} implies  that
\begin{equation}\label{instant.noW}
\dis\pa_t \CE(t,k)+\la \|\nu^{1/2}\{\FI-\FP\}\hat{\solU}\|_{L^2_{p}}^2
+\frac{\la |k|^2}{1+|k|^2}
\left(
\left|\hat{a}\right|^2
+
\left|\hat{b}\right|^2
+
\left| \hat{c} \right|^2
\right)
\\
\lesssim
 0,
\end{equation}
where note further that one has
$
\left|\hat{a}\right|^2
+
\left|\hat{b}\right|^2
+
\left| \hat{c} \right|^2
\lesssim 
\|
\FP \hat{\solU}\|_{L^2_p}^2. 
$

To do the weighted estimates, in particular for the soft potentials \eqref{hypSOFT}, we need to use the energy splitting as follows.
With \eqref{thm.tfli.p1} we define 
\begin{equation}\notag
\begin{split}
    \CE_\wN^0(t,k)\eqdef &
    \ind_{|k|\le 1} \left(
    \CE(t,k)
    +\kappa_4  \| w_{\wN}\{\FI-\FP\}\hat{\solU}\|_{L^2_p}^2 \right),
    \\
        \CE_\wN^1(t,k)\eqdef &
    \ind_{|k|> 1} \left(
    \CE(t,k)
    +\kappa_5  \| w_{\wN}\hat{\solU}\|_{L^2_p}^2 \right).
    \end{split}
\end{equation}
Again $ \kappa_4, \kappa_5>0$ will be determined later.

We prove estimates for each of these individually.  For $\CE^1_{\wN}(t,k)$ we combine \eqref{instant.noW} with \eqref{macroWeightINEQ.noM} for $|k|> 1$ to  obtain for a suitably small $\kappa_5>0$  that
$$
\partial_t  \mathcal{E}^1_{\wN}(t,k) 
+ \la  \| \nu^{1/2}w_\wN  \hat{\solU}\|_{L^2_{p}}^2 \ind_{|k|> 1} 
\lesssim
 0.
$$
Here we have used the fact that when $|k|> 1$ then 
$
\frac{ |k|^2}{1+|k|^2} \ge \frac{1}{2}.
$

Furthermore, when $|k|\le 1$ it holds that
$
\frac{ |k|^2}{1+|k|^2} \ge \frac{|k|^2}{2}.
$
In this case we combine \eqref{instant.noW} with \eqref{macroWeightINEQ}  on $|k|\le 1$ to  obtain for a small $\kappa_4>0$  that
\begin{equation}\notag
\partial_t  \mathcal{E}_{\wN}^0(t,k) 
+ \la |k|^2 \| \nu^{1/2}w_\wN  \hat{\solU}\|_{L^2_{p}}^2  \ind_{|k|\le 1} 
\lesssim
 0.
\end{equation}
Lastly we define 
$
\mathcal{E}_{\wN}(t,k) 
\eqdef
\mathcal{E}_{\wN}^0(t,k) 
+
\mathcal{E}_{\wN}^1(t,k) 
$
and we notice that \eqref{thm.tfli.1} is satisfied.
Then \eqref{thm.tfli.2} follows from adding the previous two differential inequalities.
\end{proof}

\subsection{Proof of time-decay of linear solutions}\label{sec.tf}
  Our proof of Theorem \ref{thm.ls} is based on Theorem \ref{thm.tfli} and the interpolation argument given below.

\begin{proof}[Proof of Theorem \ref{thm.ls}]   
We define 
$
  \fcnP(k)\eqdef
    \la \left( 1 \wedge |k|^2\right).
$  
By \eqref{thm.tfli.2} we have that
$$
\mathcal{E}_{\wN}(t,k)
\le
\mathcal{E}_{\wN}(0,k),
$$
for any $\wN \in \R$.
Now we use the interpolation technique as in \cite{MR2209761}, but in a different context.  In particular, for $\wE>0$, using \eqref{thm.tfli.1} we have
$$
\mathcal{E}_{\wN}(t,k)
\lesssim
 \mathcal{E}_{\wN-1}^{\wE/(\wE+1)}(t,k) ~ \mathcal{E}_{\wN+\wE}^{1/(\wE+1)}(t,k)
\lesssim
\| \nu^{1/2}w_\wN  \hat{\solU}\|_{L^2_{p}}^{2 \wE/(\wE+1)}\mathcal{E}_{\wN+\wE}^{1/(\wE+1)}(t,k).
$$
We therefore conclude that
$$
\mathcal{E}_{\wN}^{(\wE+1)/\wE}(t,k)
\lesssim
\|\nu^{1/2} w_\wN  \hat{\solU}\|_{L^2_{p}}^{2 }\mathcal{E}_{\wN+\wE}^{1/\wE}(t,k)
\lesssim
\|\nu^{1/2} w_\wN  \hat{\solU}\|_{L^2_{p}}^{2 }\mathcal{E}_{\wN+\wE}^{1/\wE}(0,k).
$$
Now we can rewrite \eqref{thm.tfli.2}, for any  $k\in \threed$, as
\begin{equation*}
    \pa_t \CE_{\wN}(t,k)+ \la \fcnP(k)  \mathcal{E}_{\wN}^{(\wE+1)/\wE}(t,k)\mathcal{E}_{\wN+\wE}^{-1/\wE}(0,k)  \leq 0.
\end{equation*}
To prove \eqref{thm.ls.1.soft}, one can bound $\CE(t,k)$ as follows
\begin{equation*}
    \pa_t \CE_{\wN}(t,k)  \mathcal{E}_{\wN}^{-1-1/\wE}(t,k)    \lesssim -\fcnP(k)  \mathcal{E}_{\wN+\wE}^{-1/\wE}(0,k).
\end{equation*}
Integrating this over time, we obtain
\begin{equation*}
    \wE \CE_{\wN}^{-1/\wE}(0,k) - \wE \mathcal{E}_{\wN}^{-1/\wE}(t,k)    \lesssim - t \fcnP(k)  \mathcal{E}_{\wN+\wE}^{-1/\wE}(0,k).
\end{equation*}
For any $\wN \in \R$ and $\wE>0$, uniformly in $k\in \threed$, we have shown that
\begin{equation*}
     \mathcal{E}_\wN(t,k)    \lesssim
    \mathcal{E}_{\wN+\wE}(0,k) 
    \left(\frac{ t \fcnP(k)}{\wE}   + 1\right)^{-\wE}.
\end{equation*}
We also just used the estimate $\mathcal{E}_{\wN}(0,k) \lesssim\mathcal{E}_{\wN+\wE}(0,k)$.

As before, we integrate over $k$ and split into $|k|\leq 1$ and $|k|> 1$ to achieve
\begin{equation*}
\int_{|k| > 1}dk~ |k|^{2m}
     \mathcal{E}_\wN(t,k)   \lesssim
     \left(\frac{ t }{\wE}   + 1\right)^{-\wE}
  \int_{|k|> 1}dk~ |k|^{2m}  \mathcal{E}_{\wN+\wE}(0,k).
\end{equation*}
Alternatively, when $|k|\leq 1$ we choose $\wE$ to be any number $\wE=\wK>2\sigma_{r,m}$ and obtain
\begin{multline*}
\int_{|k| \leq 1}dk~ |k|^{2m}
     \mathcal{E}_\wN(t,k)    
\lesssim     
     \int_{|k| \leq 1}dk~ |k|^{2m}
    \mathcal{E}_{\wN+\wK}(0,k) 
    \left(\frac{ t |k|^2}{\wK}   + 1\right)^{-\wK}
    \\
    \lesssim     
    \left( t   + 1\right)^{-2\sigma_{r,m}}
     \| w_{\wN+\wK} \solU_0\|_{L^2_pL^r_x}^2.
\end{multline*}
For $1\le r \le 2$, this uses the same H\"{o}lder and Hausdorff-Young argument.
\end{proof}

\section{Linear decay theory in $L^\infty_{p}L^2_x$}\label{I2decay}

Now we work on the linear $L^\infty_{p}L^2_x$ bounds and time decay. As is customary, we express solutions, $f(t,x,p)$, to \eqref{rBoltz} with the semigroup $U(t)$ as 
\begin{equation}
f(t,x,p) = \{U(t)f_0\}(x,p),
\label{Udef}
\end{equation}
with initial data given by
$$
\{U(0)f_0\}(x,p) = f_0(x,p).
$$
Then we have the following result.

\begin{theorem}
\label{decay0l}
Fix $\ell \geq 0,\, r\in[1,2]$ and $k\in[0,\sigma_{r,0}]$.  Suppose  
$
w_{\ell+k}f_0 \in L^\infty_{p}L^2_x,
$
then under \eqref{hypSOFT} the semi-group satisfies the estimate
$$
\|w_\ell\{U(t)f_0\} \|_{L^{\infty}_{p}L^2_x}
\le
C (1+t)^{-k}
\left( \| w_{\ell+k}f_0\|_{L^{\infty}_pL^2_x}+\|f_0\|_{L^2_pL^r_x}
\right).
$$
Above the positive constant $C= C_{\ell, k}$ only depends on $\ell$ and $k$.  
\end{theorem}

A key idea in the proof is that, instead of placing everything in the $L^\infty_{p}L^\infty_{x}$ space as in the proof of Theorem 4.1 in \cite{Strain2010}, we use the $L^\infty_{p}L^2_x$ space. The crucial new element that we now use is Minkowski's inequality in the following form
$$
\left\| \int_\Omega  f(y)  \nu(dy)  \right\|_{L^p_x}
\lesssim
 \int_\Omega  \left\|f(y)\right\|_{L^p_x}  \nu(dy),
 \quad
 1 < p < \infty,
$$
where $\Omega$ is any measure space with sigma-finite measure $\nu$.
This allows for substantial simplifications over previous work.  
Also since the inequalities used to deal with the term $H^{low,2}_5$ (which will be defined in the following discussion) in the proof of \cite[Theorem 4.1]{Strain2010} do not work well in our setting, and we use a completely different approach. The term $H^{low,2}_5$  was estimated with a  complicated change of variables in the proof of \cite[Theorem 4.1]{Strain2010}. Our new approach does not need to use any such change of variables as a result of our utilization of Minkowski's inequality.

Now we first consider solutions to the linearization of \eqref{rBoltz} with the compact operator $K$ removed from \eqref{rBoltz}. This equation is given by
\begin{gather}
\left( \partial_t  + \hat{p}\cdot \nabla_x  + \nu(p) \right) f
=
0,
\quad
f(0,x,p) =  f_0(x,p).
\label{rBlinWO}
\end{gather}
Let the semigroup 
$
G(t)f_0
$ 
denote the solution to this system \eqref{rBlinWO}.
Explicitly
$$
G(t)f_0(x,p) \eqdef e^{-\nu(p) t } f_0(x-\hat{p}t, p).
$$
For the solution $\{U(t)f_0\}(x,p)$, by iterating twice (as did Vidav \cite{MR0259662}) we have
\begin{gather}
\notag
\left\{  U(t) f_0 \right\}(x,p) =
G(t)f_0(x,p)
+
\int_0^t  ~ ds_1 ~ G(t-s_1)K^{1-\chi} \left\{  U(s_1) f_0 \right\}(x,p)
\\
+
\int_0^t  ~ ds_1 ~G(t-s_1)K^{\chi}  G(s_1) f_0 (x,p)
\label{linearEXP}
\\
\notag
+
\int_0^t ~ ds_1 ~  
\int_0^{s_1} ~ ds_2 ~
G(t-s_1)K^{\chi}
G(s_1-s_2)
K^{1-\chi}
\left\{  U(s_2) f_0 \right\}(x,p)
\\
+
\int_0^t ~ ds_1 ~  
\int_0^{s_1} ~ ds_2 ~
G(t-s_1)K^{\chi}
G(s_1-s_2)
K^{\chi}
\left\{  U(s_2) f_0 \right\}
(x,p).
\notag
\end{gather}
Equivalently
\begin{equation}
\left\{  U(t) f_0 \right\}(x,p) 
\eqdef
H_1(t,x,p)
+
H_2(t,x,p)
+
H_3(t,x,p)
+
H_4(t,x,p)
+
H_5(t,x,p),
\notag
\end{equation}
where
\begin{eqnarray*}
H_1(t,x,p)
& \eqdef &
e^{-\nu(p) t } f_0(x-\hat{p}t, p),
\\
H_2(t,x,p)
& \eqdef &
\int_0^t  ~ds_1 ~ e^{-\nu(p) (t-s_1) }
K^{1-\chi} \left\{  U(s_1) f_0 \right\}(y_1,p),
\\
H_3(t,x,p)
& \eqdef &
\int_0^t  ~ds_1 ~ e^{-\nu(p) (t-s_1) }
\int_{\mathbb{R}^3} dq_1~ k^\chi(p,q_1)~ e^{-\nu(q_1) s_1 }
f_0(y_1-\hat{q}_1s_1,q_1).
\end{eqnarray*}
Just above and below we will be using the following short hand notation
\begin{eqnarray}
\notag
y_1 & \eqdef & x- \hat{p}(t-s_1),
\\
y_2 & \eqdef & y_1-\hat{q}_1(s_1-s_2)
=
x- \hat{p}(t-s_1)-\hat{q}_1(s_1-s_2).
\label{zDEF}
\end{eqnarray}
We are also using the notation
$
q_1^0 = \sqrt{1+|q_{1}|^2}
$
and
$
q_1= (q_{11},q_{12},q_{13})\in\mathbb{R}^3
$
with
$
\hat{q}_1 = q_1/q_1^0.
$
Furthermore the next term is
\begin{multline*}
H_4
(t,x,p)
\eqdef
\int_{\mathbb{R}^3} dq_1~ k^\chi(p,q_1)
\int_0^t  ~ds_1 \int_0^{s_1}  ~ds_2 ~
e^{-\nu(p) (t-s_1) }
 e^{-\nu(q_1) (s_1- s_2) }
 \\
 \times
K^{1-\chi}
\left\{  U(s_2) f_0 \right\}
(y_2,q_1).
\end{multline*}
Lastly, we may also expand out the fifth component as
\begin{multline}
H_5(t,x,p)
=
\int_{\mathbb{R}^3} dq_1~ k^\chi(p,q_1)~ 
\int_{\mathbb{R}^3} dq_2~ k^\chi(q_1,q_2)~ 
\int_0^t  ~ds_1 ~ e^{-\nu(p) (t-s_1) }
\\
\times
\int_0^{s_1} ~ ds_2 ~
e^{-\nu(q_1) (s_1-s_2) }
\left\{  U(s_2) f_0 \right\}(y_2,q_2).
\label{h5}
\end{multline}

We have the following estimates for the $H_i$ terms above ($1\leq i\leq 5$).

\begin{lemma}
\label{hESTIMATES}
Given $\ell \ge 0$,
for any  $k\ge 0$ we have
\begin{equation*}
\|
w_\ell
H_1
\|_{L^\infty_{p}L^2_x}
+
\|
w_\ell
H_3
\|_{L^\infty_{p}L^2_x}
\le C_{\ell,k} (1+t)^{-k}\|w_{k+\ell} f_0\|_{L^\infty_pL^2_x}.
\end{equation*}
\end{lemma}

\begin{lemma}
\label{hESTIMATE24}
Fix  $\ell \ge 0$.
For any small $\eta>0$ and any $k\ge 0$  we have
\begin{equation*}
\|w_\ell
H_2\|_{L^\infty_{p}L^2_x}
+
\|
w_\ell
H_4
\|_{L^\infty_{p}L^2_x}
\le \eta(1+t)^{-k}||\varpi _{k} w_\ell f ||_{L^\infty_{p,t}L^2_x}.
\end{equation*}
\end{lemma}

\begin{lemma}
\label{hESTIMATE3}
Fix $\ell \ge 0$, choose any (possibly large) $j>0$.
For any small $\eta>0$ and any $k\ge 0$ we have the estimate
\begin{multline*}
\|
w_\ell
H_5
\|_{L^\infty_{p}L^2_x}
\le \eta(1+t)^{-k} ||\varpi_k w_\ell f ||_{L^\infty_{p,t}L^2_x}\\
+
C_\eta 
\int_0^t ds ~ e^{-\eta (t-s)} \|w_{-j}f\|_{L^2_{x,p}}(s)
+
\|w_\ell
R_1(f)(t)
\|_{L^\infty_{p}L^2_x}.
\notag
\end{multline*}
By the $L^2_p L^2_x$ decay theory from Theorem \ref{thm.ls} and Proposition \ref{BasicDecay} (below), we can choose $j>\max \{2\sigma_{r,0},6/b+k\}$ and take $k\in [0,\sigma_{r,0}]$ to have   
\begin{multline*}
\int_0^t ds ~ e^{-\eta (t-s)} \|w_{-j}f\|_{L^2_{p}L^2_{x}}(s)\\
\le
C_\eta (1+t)^{-k}\|w_{2k-j}f_0 \|_{L^2_{p}L^2_{x}}+C_\eta (1+t)^{-k}\|f_0 \|_{L^2_pL^r_x}\\
\le
C_\eta (1+t)^{-k}
\left( \|w_{\ell+k}f_0 \|_{L^\infty_pL^2_x}+\|f_0 \|_{L^2_pL^r_x} \right).
\end{multline*}
The above estimates hold for any $k\in [0,\sigma_{r,0}]$ and $\ell\geq 0$.  On the other hand, for the last term
 involving $R_1(f)$ if we restrict $k\in[0,1]$ then $\forall \eta>0$ we have
$$
\|w_\ell
R_1(f)
\|_{L^\infty_{p}L^2_x}
\le
\eta(1+t)^{-k}||\varpi_k w_\ell f ||_{L^\infty_{p,t}L^2_x}.
$$
This term $R_1$ is defined in (\ref{r1def}) during the course of proof.
\end{lemma}

Combining Lemma \ref{hESTIMATES}, Lemma \ref{hESTIMATE24} and Lemma \ref{hESTIMATE3} gives Theorem \ref{decay0l}. Actually, Lemma \ref{hESTIMATES}, Lemma \ref{hESTIMATE24} and Lemma \ref{hESTIMATE3} imply that
for any $\eta >0$ and $k\in [0,\sigma_{r,0}]$ we have
\begin{multline*}
\|w_\ell f \|_{L^\infty_p L^2_x}(t)
\le
C_{\ell,k,\eta} (1+t)^{-k}\| w_{\ell+k}f_0\|_{L^\infty_p L^2_x}\\
+
\eta (1+t)^{-k} ||\varpi_k w^\ell f ||_{L^\infty_{p,t} L^2_x}
+
C_\eta (1+t)^{-k}\|f_0 \|_{L^2_pL^r_x}.
\end{multline*}
Equivalently
$$
 ||\varpi_k w^\ell f ||_{L^\infty_{p,t} L^2_x}
\le
C_{\ell,k,\eta}\| w_{\ell+k}f_0\|_{L^\infty_p L^2_x}
+
C_\eta \|f_0 \|_{L^2_pL^r_x}.
$$
With this inequality, we have proved Theorem \ref{decay0l} subject to  Lemmas \ref{hESTIMATES}, \ref{hESTIMATE24}, and \ref{hESTIMATE3}.  We now prove those lemmas. \\

We will use the following known basic decay estimate, as in \cite[Proposition 4.5]{Strain2010}:
\begin{proposition}
\label{BasicDecay}
Suppose without loss of generality that $\lambda \ge \mu \ge 0$.  Then
$$
\int_0^t  \frac{ds}{(1+ t - s)^{\lambda} (1+ s)^{\mu}}
\le
\frac{C_{\lambda,\mu}(t)}{ (1+t)^{\rho}},
$$
where $\rho = \rho(\lambda,\mu) =\min\{ \lambda+\mu -1, \mu\}$ and
$$
0 \le 
C_{\lambda,\mu}(t) = 
C
\left\{
\begin{array}{cc}
1 & \mbox{if} ~ \lambda \ne 1,
\\
\log(2+t) & \mbox{if} ~ \lambda = 1.
\end{array}
\right.
$$
\end{proposition}

Furthermore, we will use the following basic estimate from the Calculus
\begin{equation}
e^{-ay} (1+y)^k \le \max\{1, e^{a-k} k^k a^{-k}\},
\quad
a,y, k \ge 0.
\label{ElemCalc}
\end{equation}
We are now ready to prove the lemmas above.

\begin{proof}[Proof of Lemma \ref{hESTIMATES}.]
We first look at $H_1$. The following result is shown in the proof of Lemma 4.2 in \cite{Strain2010}:
\begin{equation}
e^{-\nu(p) t }
\le
C_{k} 
p_0^{kb/2}
(1+t)^{-k}
\le
C_{k} 
w_k(p)
(1+t)^{-k},
\quad
\forall t, k >0.
\label{polyE}
\end{equation}
So we have 
\begin{multline*}
\left|
w_\ell(p) H_1(t,x,p)
\right|
=
\left|
w_\ell(p) ~ e^{-\nu(p) t } f_0(x-\hat{p}t, p)
\right|\\
\le
C_{k}\left|
w_\ell(p)w_k(p)
(1+t)^{-k} f_0(x-\hat{p}t, p)
\right|.
\end{multline*}
So
$$
\|
w_\ell H_1
\|_{L^2_x}
\le
C_{k}\|
w_\ell w_k
\varpi_{-k} f_0
\|_{L^2_x}.
$$
Thus
 \begin{equation*}
\|
w_\ell
H_1
\|_{L^\infty_{p}L^2_x}
\le C_{\ell,k} (1+t)^{-k}\|w_{k+\ell} f_0\|_{L^\infty_pL^2_x}.
\end{equation*}
Now we turn to $H_3$. Since 
$$
 e^{-\nu(p) (t-s_1) }
  e^{-\nu(q) s_1 }
  \le
  e^{-\nu(\max\{|p|,|q|\}) t },
  $$
  where 
  $
  \nu(\max\{|p|,|q|\})
$
is $\nu$ evaluated at 
$
\max\{|p|,|q|\},
$
we have
\begin{multline*}
\left| 
w_\ell (p)
H_3(t,x,p)
\right|\\
\le
w_\ell (p)
\int_0^t  ~ds_1 ~ e^{-\nu(p) t }
\int_{|p| \ge |q_1|} dq_1 \left| k^\chi(p,q_1) \right| 
\left| 
f_0(x- \hat{p}(t-s_1)-\hat{q}_1 s_1,q_1)
\right|
\\
+
w_\ell (p)
\int_0^t  ~ds_1 ~ 
\int_{|p| < |q_1|} dq_1~ \left| k^\chi(p,q_1) \right| ~ e^{-\nu(q_1) t }
\left| 
f_0(x- \hat{p}(t-s_1)-\hat{q}_1 s_1,q_1)
\right|.
\end{multline*}
So by Minkowski's inequality, we have
\begin{multline*}
\| 
w_\ell 
H_3
\|_{L^2_x}
\le
w_\ell (p)
\int_0^t  ~ds_1 ~ e^{-\nu(p) t }
\int_{|p| \ge |q_1|} dq_1 \left| k^\chi(p,q_1) \right| 
\| 
f_0
\|_{L^2_x}(q_1)
\\
+
w_\ell (p)
\int_0^t  ~ds_1 ~ 
\int_{|p| < |q_1|} dq_1~ \left| k^\chi(p,q_1) \right| ~ e^{-\nu(q_1) t }
\| 
f_0
\|_{L^2_x}(q_1).
\end{multline*}
We will estimate the second term, and we  remark that the first term can be handled in exactly the same way.
As in Lemma 4.2 in \cite{Strain2010}, and similar to \eqref{polyE},  we have
$$
e^{-\nu(q_1) t }
\| 
f_0
\|_{L^2_x}(q_1)\le
C_k(1+t)^{-k-1} w_{k+1}(q_1)\| 
f_0
\|_{L^2_x}(q_1).
$$
Next we use the estimate for $k^\chi(p,q)$ from Lemma \ref{boundK2}.
When using this estimate
 we may suppose
$
|q_1| \le 2|p|.
$
For otherwise, if say $|q_1| \ge 2 |p|$, then we have 
$
|p-q_1| \ge |q_1| /2
$
which leads directly to
\begin{gather}
w_{k+1}(q_1) w_\ell (p)
e^{-c|p-q_1|/2}
\le
Cw_{k+1+\ell}(q_1) 
e^{-c|q_1|/4}
\le C.
\label{comparablePQ}
\end{gather}
In this case we easily obtain the following estimate:
\begin{multline*}
\int_0^t  ~ds_1 ~ 
w_\ell (p)
\int_{|p| < |q_1|} dq_1~{\bf 1}_{|q_1| \ge 2|p|}~ \left| k^\chi(p,q_1) \right| ~ e^{-\nu(q_1) t }
\| 
f_0
\|_{L^2_x}(q_1)
\\
\lesssim
(1+t)^{-k} \| w_{-j}f_0 \|_{L^\infty_p L^2_x},
\quad \forall k>0, ~ j >0.
\end{multline*}
Thus in the following  we assume 
$
 |p| < |q_1| \le 2 |p|.
$
On this region we have
\begin{multline}
\int_0^t  ~ds_1 ~ 
w_\ell (p)
\int_{|p| < |q_1|\le 2|p|} dq_1~ \left| k^\chi(p,q_1) \right|  ~ e^{-\nu(q_1) t }
\| 
f_0
\|_{L^2_x}(q_1)
\\
\le
C_{k} \frac{t}{(1+t)^{k+1}}
\int_{|p| < |q_1|\le 2|p|} dq_1~
w_\ell (q_1) ~ w_{k+1 }(q_1) ~ \left| k^\chi(p,q_1) \right| ~
\| 
f_0
\|_{L^2_x}(q_1)
\\
\le
C_{k}(1+t)^{-k} \| w_{\ell+k}f_0\|_{L^\infty_p L^2_x}
~ w_{1 }(p)
\int_{|p| < |q_1|\le 2|p|} dq_1~ \left| k^\chi(p,q_1) \right|.
\label{desired1}
\end{multline}
Then from Lemma \ref{boundK2} we clearly have the following bound 
$$
w_{1 }(p)
\int_{|p| < |q_1|\le 2|p|} dq_1~ \left| k^\chi(p,q_1) \right|
\le C.
$$
This completes the time decay estimate for $H_3$ and our proof of  the lemma.
\end{proof}

\noindent {\it Proof of Lemma \ref{hESTIMATE24}.}
As in the proof of Lemma 4.3 in \cite{Strain2010}, we will use the following lemma (Lemma 4.6 in \cite{Strain2010}).
\begin{lemma}[\cite{Strain2010}]
\label{boundKinfX}  
Fix any $\ell \ge 0$ and any $j>0$.
Then given any small $\eta >0$, which depends upon $\chi$ in \eqref{cut}, the following estimate holds
$$
\left| w_\ell(p) K_i^{1-\chi}(h)(p) \right| 
\le
\eta e^{- cp^0} \| w_{-j}h\|_{L^\infty_{q}}.
$$
Above the constant $c>0$ is independent of $\eta$ and $i=1,2$.
\end{lemma}

From \eqref{kCUT} and Minkowski's inequality, we have
$$\|K_i^{1-\chi}(h)(p)\|_{L^2_x}\le K_i^{1-\chi}(\|h\|_{L^2_x})(p).$$
So with the lemma above, we have 
\begin{equation}\label{newlemma4.6}
\| w_\ell(p) K_i^{1-\chi}(h)(p) \|_{L^2_x} 
\le
\eta e^{- cp^0} \|w_{-j} h\|_{L^\infty_q L^2_x}.
\end{equation} 
Notice that in the result above, we use $ \|w_{-j} h\|_{L^\infty_q L^2_x}$ instead of $\|w_{-j} h\|_{L^\infty_p L^2_x}$ to avoid confusion, since $p$ is a variable in the left-hand-side. From now on, we always use $q$ to denote the momentum variable in the norm expression whenever $p$ has been used to denote a particular momentum value in an inequality.
For $H_2$, with the result above, for any small $\eta'>0$ we have
\begin{multline*}
\|
w_\ell H_2 
\|_{L^2_x}(t,p)
=
w_\ell(p)
\left\| 
\int_0^t  ~ds_1 ~ e^{-\nu(p) (t-s_1) }
K^{1-\chi} \left\{  U(s_1) f_0 \right\}(y_1,p)
\right\|_{L^2_x}
 \\
\le
\eta'  ~ e^{-c p^0} || \varpi_k  f||_{L^\infty_{q,t}L^2_x}
 \int_0^t  ~ds_1 ~ e^{-\nu(p) (t-s_1) } (1+s_1)^{-k}, \quad \forall j>0.
\end{multline*}
For any $\lambda >\max\{1,k\}$ 
we have
\begin{multline*}
\|
w_\ell H_2 
\|_{L^2_x}(t,p)\\
\le
\eta'  ~ w_{\lambda}(p) e^{-c p^0} || \varpi_k  f||_{L^\infty_{q,t}L^2_x}
 \int_0^t  ~ds ~ (1+t-s)^{-\lambda} (1+s)^{-k}
\\
\le
\eta ~ (1+t)^{-k}
|| \varpi_k  f||_{L^\infty_{q,t}L^2_x},
\end{multline*}
which follows from Proposition \ref{BasicDecay}.
This is the desired estimate for $H_2$.
For $H_4$ we once again use \eqref{newlemma4.6}, for any small $\eta'>0$, to obtain
\begin{multline*}
\|
w_\ell
H_4
\|_{L^2_x}(t,p)
\le
\int_{\mathbb{R}^3} dq_1~ k^\chi(p,q_1)
\int_0^t  ~ds_1 \int_0^{s_1}  ~ds_2 ~
e^{-\nu(p) (t-s_1) }
\\
\times
 e^{-\nu(q_1) (s_1- s_2) }
w_\ell(p)
\|
K^{1-\chi}\left(
\left\{  U(s_2) f_0 \right\}
\right)
(y_2,q_1)
\|_{L^2_x}
\\
\le
\eta' || \varpi_k  f ||_{L^\infty_{q,t}L^2_x}
\int_{\mathbb{R}^3} dq_1~ k^\chi(p,q_1)
\frac{w_\ell(p)}{w_\ell(q_1)}
e^{-cq_1^0}
\\
\times
 \int_0^t  ~ds_1 ~ 
\int_0^{s_1}  ~d s_2 ~ e^{-\nu(p) (t-s_1) } e^{-\nu(q_1) (s_1- s_2) }(1+s_2)^{-k}.
\end{multline*}
For the time decay, from Proposition \ref{BasicDecay}, 
we notice that
\begin{multline}
 \int_0^t  ~ds_1 ~ 
\int_0^{s_1}  ~d s_2 ~ e^{-\nu(p) (t-s_1) } e^{-\nu(q_1) (s_1- s_2) }(1+s_2)^{-k}
\\
\lesssim
w_{ \lambda}(p) w_{\lambda}(q_1)
\int_0^t  ~ds_1 ~ (1+t-s_1)^{-\lambda}
\int_0^{s_1}  ~ds_2 ~ (1+s_1-s_2)^{-\lambda}(1+s_2)^{-k}
\\
\lesssim
 w_{\lambda}(p) w_{\lambda}(q_1)
(1+t)^{-k}.
\notag
\end{multline}
Above we have taken $\lambda >\max\{1,k\}$.
Combining these estimates yields
\begin{align*}
&\quad \|
w_\ell(p)
H_4 
\|_{L^2_x}\\
&\lesssim
 \eta' (1+t)^{-k} ||\varpi_k w_\ell f ||_{L^\infty_{q,t}L^2_x} 
\int_{\mathbb{R}^3} dq_1~ k^\chi(p,q_1)
w_{\ell + \lambda}(p) w_{-\ell + \lambda}(q_1)
e^{-cq_1^0}.
\end{align*}
To estimate the remaining integral and weights
we split into three cases.
If either
$
2|q_1| \le |p|,
$
or
$
|q_1| \ge 2|p|,
$
then we bound all the weights and the remaining momentum integral by a constant as in 
\eqref{comparablePQ}.
Alternatively if
$
\frac{1}{2} |q_1|
\le
|p| 
\le 
2|q_1|,
$
then the desired estimate is obvious since we have strong exponential decay in both $p$ and $q_1$.  
In either of these  cases we have the estimate for $H_4$. 
\qed \\

\noindent {\it Proof of Lemma \ref{hESTIMATE3}.}
We will utilize rather extensively the estimate for $k^\chi$ from Lemma \ref{boundK2}.
We now further split $H_5(t,x,p)$ as
\begin{equation}
H_5(t,x,p) =  H_5^{high}(t,x,p)+ H_5^{low}(t,x,p),
\label{h5sum}
\end{equation}
and estimate each term on the right individually.    For $M>>1$ 
we define
\begin{equation}
{\bf 1}_{high} 
\eqdef
{\bf 1}_{|p|> M}  {\bf 1}_{|q_1|\le M} 
+ 
 {\bf 1}_{|q_1|> M}.
 \label{1high}
\end{equation}
Notice
$
{\bf 1}_{high} +{\bf 1}_{|p|\le M}  {\bf 1}_{|q_1|\le M} 
=1.
$
Now the first term in the expansion is 
\begin{gather*}
H_5^{high}(t,x,p)
\eqdef
\int_{\mathbb{R}^3}  dq_1 ~ k^\chi(p,q_1)~ 
\int_{\mathbb{R}^3} dq_2~ k^\chi(q_1,q_2)~  {\bf 1}_{high} ~
\int_0^t  ~ds_1 ~ e^{-\nu(p) (t-s_1) }
\\
\times
\int_0^{s_1} ~ ds_2 ~
e^{-\nu(q_1) (s_1- s_2) }
\left\{  U(s_2) f_0 \right\}(y_2,q_2).
\end{gather*}
The proof of Lemma 4.4 in \cite{Strain2010}, and also \eqref{polyE}, shows that for any $\lambda \ge 0$ we have
\begin{multline*}
\int_0^t  ~ds_1 ~ e^{-\nu(p) (t-s_1) }
\int_0^{s_1} ~ ds_2 ~
e^{-\nu(q_1) (s_1- s_2) }
\\
\le
C_\lambda w_\lambda (p) w_\lambda (q_1)
\int_0^t  ~ds_1 ~ 
\int_0^{s_1} ~ ds_2 ~
(1+(t-s_1))^{-\lambda} (1+(s_1 - s_2))^{-\lambda}.
\end{multline*}
When either $|p|> M$ or $|q_1|> M$, by Lemma \ref{boundK2}, we have the bound
$$
\left| k^{\chi}(p,q_1)  \right|
\le C M^{-\zeta} 
\left( p^0+ q_1^0 \right)^{-b/2}
e^{-c |p-q_1|}.
$$
If either $|p| \ge 2|q_1|$ or $|q_1| \ge 2|p|$
 then as in \eqref{comparablePQ} we have
$$
w_{\ell+\lambda}(p) w_\lambda(q_1) e^{-c|p-q_1|} 
\le C.
$$
Thus by combining the last few estimates we have
\begin{multline*}
w_\ell(p) \int_{\mathbb{R}^3}  dq_1 ~ \left| k^{\chi}(p,q_1)  \right|~ 
\int_{\mathbb{R}^3} dq_2~  \left| k^{\chi}(q_1,q_2)  \right|~ {\bf 1}_{high} ~
\left( {\bf 1}_{|p| \ge 2|q_1|}+{\bf 1}_{|p| \le \frac{1}{2}|q_1|} \right)
\\
\times
\int_0^t  ~ds_1 ~ e^{-\nu(p) (t-s_1) }
\int_0^{s_1} ~ ds_2 ~
e^{-\nu(q_1) (s_1- s_2) }
\|
\left\{  U(s_2) f_0 \right\}(y_2,q_2)
\|_{L^2_x}
\\
\le
\frac{C_{\lambda}}{M^{\zeta+b/2}} \| \varpi_k  w_\ell f \|_{L^\infty_{q,t}L^2_x} \int_{\mathbb{R}^3}  dq_1 ~ e^{-c|p-q_1|}~ 
\int_{\mathbb{R}^3} dq_2~ e^{-c|q_1-q_2|}~ 
\\
\times
\int_0^t  ~ds_1 ~ 
\int_0^{s_1} ~ ds_2 ~
\frac{(1+s_2)^{-k}}{(1+(t-s_1))^{\lambda} (1+(s_1 - s_2))^{\lambda} }.
\end{multline*}
With Proposition \ref{BasicDecay}, for any $k\ge 0$ and $\lambda> \max\{k,1\}$ the previous term is bounded from above by
  \begin{gather*}
\le
\frac{C_{k,\lambda}}{M^{\zeta+b/2}}~ (1+t)^{-k} ~ || \varpi_k w_\ell f||_{L^\infty_{q,t}L^2_x}.
\end{gather*}
This is the desired estimate for $M\gg 1$ chosen sufficiently large.
We now consider the remaining part of $H_5^{high}$.
As in the previous estimates and \eqref{comparablePQ}, if either 
 $|q_2 | \ge 2|q_1|$ or
 $|q_1| \ge 2|q_2|$ then for any $k\ge 0$ we have
 \begin{multline*}
w_\ell(p) \int_{\mathbb{R}^3}  dq_1 ~ \left| k^{\chi}(p,q_1)  \right|~ 
\int_{\mathbb{R}^3} dq_2~  \left| k^{\chi}(q_1,q_2)  \right|~
 \left( {\bf 1}_{|q_2| \ge 2|q_1|}+{\bf 1}_{|q_2| \le \frac{1}{2}|q_1|} \right)
\\
\times
{\bf 1}_{\frac{1}{2}|p| \le |q_1|\le 2|p|}
{\bf 1}_{high} 
\int_0^t  ~ds_1 ~ e^{-\nu(p) (t-s_1) }
\int_0^{s_1} ~ ds_2 ~
e^{-\nu(q_1) (s_1- s_2) }
\\
\times
\| 
\left\{  U(s_2) f_0 \right\}(q_2)
\|_{L^2_x}\\
\le
\frac{C_{k,\lambda}}{M^{2\zeta+b}}|| \varpi_k w_\ell f||_{L^\infty_{q,t}L^2_x} \int_{\mathbb{R}^3}  dq_1 ~ e^{-c|p-q_1|}~ 
\int_{\mathbb{R}^3} dq_2~ e^{-c|q_1-q_2|}~ 
\\
\times
\int_0^t  ~ds_1 ~ 
\int_0^{s_1} ~ ds_2 ~
\frac{(1+s_2)^{-k}}{(1+(t-s_1))^{\lambda} (1+(s_1 - s_2))^{\lambda} }
\\
\le
\frac{C_{k,\lambda}}{M^{2\zeta+b}}~ (1+t)^{-k} ~ || \varpi_k w_\ell f||_{L^\infty_{q,t}L^2_x}.
\end{multline*}
Above we have used exactly the same estimates as in the prior case.  Both of the last two terms have a suitably small constant in front if $M$ is  sufficiently large.
Thus the remaining part of $H_5^{high}$ to estimate is $R_1(f)(t)$ which is defined by 
\begin{gather}
\notag
R_1(f)(t)\eqdef \int_{\mathbb{R}^3}  dq_1 ~ k^\chi(p,q_1)~ 
\int_{\mathbb{R}^3} dq_2~ k^\chi(q_1,q_2)~ 
{\bf 1}_{\frac{1}{2}|p| \le |q_1|\le 2|p|} {\bf 1}_{\frac{1}{2}|q_1| \le |q_2|\le 2|q_1|} 
\\
\times
{\bf 1}_{high} ~
\int_0^t  ~ds_1 ~ e^{-\nu(p) (t-s_1) }
\int_{0}^{s_1} ~ ds_2 ~
e^{-\nu(q_1) (s_1- s_2) }
\left\{  U(s_2) f_0 \right\}(y_2,q_2).
\label{r1def}
\end{gather}
 Since all the momentum variables are comparable, we have
\begin{multline}
\notag
\| R_1(f)(t) \|_{L^2_x} 
\le 
\int_{\mathbb{R}^3}  dq_1 ~ \left| k^{\chi}(p,q_1)  \right|~ 
\int_{\mathbb{R}^3} dq_2~  \left| k^{\chi}(q_1,q_2)  \right|~
\\
\times
{\bf 1}_{\frac{1}{2}|p| \le |q_1|\le 2|p|} {\bf 1}_{\frac{1}{2}|q_1| \le |q_2|\le 2|q_1|}
{\bf 1}_{high}
\\
\times
\int_0^t  ~ds_1 ~ e^{-c\nu(q_1) (t-s_1) }
\int_0^{s_1} ~ ds_2 ~
e^{-c\nu(q_1) (s_1-s_2) }
\| 
\left\{  U(s_2) f_0 \right\}(y_2,q_2)
\|_{L^2_x}.
\notag
\end{multline}
Next using similar techniques as in the previous two estimates, including \eqref{polyE}, we obtain the following upper bound for any $k\in[0,1]$:
\begin{multline*}
w_\ell(p)\| R_1(f)(t) \|_{L^2_x} 
\le 
\frac{C_{k,\lambda}}{M^{2\zeta}}|| \varpi_k w_\ell f||_{L^\infty_{q,t}L^2_x}  \int  dq_1 ~ (q_1^0)^{-b/2 - \zeta}e^{-c|p-q_1|}~ 
\\
\times
\int dq_2~  (q_1^0)^{-b/2 - \zeta}e^{-c|q_1-q_2|}~
w_{2+2\delta}(q_1)
\int_0^{t} ~ \frac{ ds_1 }{ (1+(t - s_1))^{1+\delta} }
\\
\times
\int_0^{s_1} ~ 
\frac{ ds_2 }{ (1+(s_1 - s_2))^{1+\delta} (1+s_2)^{k}}\\
\le 
\frac{C_{k,\lambda}}{M^{2\zeta}}
(1+t)^{-k}
|| \varpi_k w_\ell f||_{L^\infty_{q,t}L^2_x}.
\end{multline*}
In this computation we have chosen $\delta >0$ to satisfy $\delta < \zeta$ where $\zeta>0$ is defined in the statement of Lemma \ref{boundK2}.  This choice guarantees that
$
w_{2+2\delta}(q_1)(q_1^0)^{-b - 2 \zeta} \le C.
$

We are ready to define the second term in our splitting of $H_5$.  It must be 
\begin{multline*}
H_5^{low}(t,x,p)
\eqdef
{\bf 1}_{|p|\le M} 
\int_{|q_1|\le M}  dq_1 ~ k^\chi(p,q_1)~ 
\int_{\mathbb{R}^3} dq_2~ k^\chi(q_1,q_2)~ 
\\
\times
\int_0^t  ~ds_1 ~ e^{-\nu(p) (t-s_1) }
\int_0^{s_1} ~ ds_2 ~
e^{-\nu(q_1) (s_1- s_2) }
\left\{  U(s_2) f_0 \right\}(y_2,q_2).
\end{multline*}
For any small $\kappa>0$, we further split this term into two terms, one of which is
\begin{multline*}
H_5^{low,\kappa}(t,x,p)
\eqdef
\int_{|q_1|\le M}  dq_1 ~ k^\chi(p,q_1)~ 
\int_{\mathbb{R}^3} dq_2~ k^\chi(q_1,q_2)~ 
\int_0^{\kappa}  ~ds_1 ~ e^{-\nu(p) (t-s_1) }
\\
\times
{\bf 1}_{|p|\le M} 
\int_0^{s_1} ~ ds_2 ~
e^{-\nu(q_1) (s_1- s_2) }
\left\{  U(s_2) f_0 \right\}(y_2,q_2)
\\
+
\int_{|q_1|\le M}  dq_1 ~ k^\chi(p,q_1)~ 
\int_{\mathbb{R}^3} dq_2~ k^\chi(q_1,q_2)~ 
\int_{\kappa}^t  ~ds_1 ~ e^{-\nu(p) (t-s_1) }
\\
\times
{\bf 1}_{|p|\le M} 
\int_{s_1-\kappa}^{s_1} ~ ds_2 ~
e^{-\nu(q_1) (s_1- s_2) }
\left\{  U(s_2) f_0 \right\}(y_2,q_2).
\end{multline*}
The other term in this latest splitting is defined just below as $H_5^{low,2}$.
Since $p$ and $q_1$ are both bounded by $M$, from Lemma \ref{nuEST} 
we have
\begin{gather}
{\bf 1}_{|p|\le M}  {\bf 1}_{|q_1|\le M}  ~
 e^{-\nu(p) (t-s_1) -\nu(q_1) (s_1- s_2) }
\le
e^{-C(t- s_2)/M^{b/2} }.
\label{MdecayEST}
\end{gather}
Then for the first term in $H_5^{low,\kappa}$ above multiplied by $w_\ell(p)$  we have the bound
\begin{align*}
&\quad w_\ell(p)\int_{|q_1|\le M} dq_1 ~ \left| k^{\chi}(p,q_1)  \right|~ 
\int_{\mathbb{R}^3} dq_2~  \left| k^{\chi}(q_1,q_2)  \right|~
\int_0^{\kappa}  ~ds_1 ~ e^{-\nu(p) (t-s_1) }
\\
&\quad\times
{\bf 1}_{|p|\le M} 
\int_0^{s_1} ~ ds_2 ~
e^{-\nu(q_1) (s_1- s_2) }
\|\left\{  U(s_2) f_0 \right\}(q_2)\|_{L^2_x}
\\
&\le
C_M || w_\ell f ||_{L^\infty_{q,t}L^2_x} 
\int_0^{\kappa}  ~ds_1 ~\int_0^{s_1} ~ ds_2 ~
 e^{-C(t- s_2)/M^{b/2} }
 \\
&\le
C_M\kappa^2  || w_\ell f ||_{L^\infty_{q,t}L^2_x}  ~ e^{-Ct/M^{b/2} } e^{C\kappa /M^{b/2} }
 \\
&\le
C_M \kappa^2  (1+t)^{-k}  || \varpi_k w_\ell f  ||_{L^\infty_{q,t}L^2_x}.
\end{align*}
 We obtain the desired estimate for the above terms by first choosing $M$ large, and second choosing $\kappa=\kappa(M)>0$ sufficiently small.

For the second term in $H_5^{low,\kappa}$ multiplied by $w_\ell(p)$  for any $k\ge 0$ we have
\begin{align*}
&\quad w_\ell(p)
\int_{|q_1|\le M} dq_1 ~ \left| k^{\chi}(p,q_1)  \right|~ 
\int_{\mathbb{R}^3} dq_2~  \left| k^{\chi}(q_1,q_2)  \right|~
\int_{\kappa}^t  ~ds_1 ~ e^{-\nu(p) (t-s_1) }
\\
&\quad\times
{\bf 1}_{|p|\le M} 
\int_{s_1-\kappa}^{s_1} ~ ds_2 ~
e^{-\nu(q_1) (s_1- s_2) }
\|\left\{  U(s_2) f_0 \right\}(q_2)\|_{L^2_x}
\\
&\le
C_M ||\varpi_k w_\ell f ||_{L^\infty_{q,t}L^2_x}\\
&\quad\times
\int_{\kappa}^t  ~ds_1 ~ \int_{s_1-\kappa}^{s_1} ~ ds_2 ~
 e^{-C(t- s_1)/M^{b/2} }e^{-C(s_1- s_2)/M^{b/2} }(1+s_2)^{-k}.
\end{align*}
Since $s_2\in [s_1 - \kappa, s_1]$ and $\kappa\in(0,1/2)$, we know
$
(1+s_2) \ge 
\left(\frac{1}{2}+s_1\right).
$
Then the previous display is further bounded above as
\begin{align*}
&\le
C_M \kappa ||\varpi_k w_\ell f ||_{L^\infty_{q,t}L^2_x}
\int_{\kappa}^t  ~ds_1 ~ 
 e^{-C(t- s_1)/M^{b/2} }(1+s_1)^{-k}
 \\
 &\le
C_M \kappa (1+t)^{-k} ||\varpi_k w_\ell f ||_{L^\infty_{q,t}L^2_x}.
\end{align*}
In the last step we have used Proposition \ref{BasicDecay}.  
We conclude the desired estimate for $H_5^{low,\kappa}$ by first choosing $M$ large, and  then $\kappa>0$ sufficiently small.

The only remaining part of $
H_5^{low}(t,x,p)
$
to be estimated  is given by 
\begin{multline*}
H_5^{low,2}(t,x,p)
\eqdef
\int_{|q_1|\le M}  dq_1 ~ k^\chi(p,q_1)~ 
\int_{\mathbb{R}^3} dq_2~ k^\chi(q_1,q_2)~ 
\int_{\kappa}^t  ~ds_1 ~ e^{-\nu(p) (t-s_1) }
\\
\times
{\bf 1}_{|p|\le M} 
\int_0^{s_1-\kappa} ~ ds_2 ~
e^{-\nu(q_1) (s_1- s_2) }
\left\{  U(s_2) f_0 \right\}(y_2,q_2),
\end{multline*}
Now for $H_5^{low,2}(t,x,p)$, we notice using Cauchy-Schwartz that
\begin{multline*}
\left|
\int_{|q_1|\le M}  dq_1 ~ k^\chi(p,q_1)~ 
\int_{\mathbb{R}^3} dq_2~ k^\chi(q_1,q_2)~ 
\|\{U(s_2) f_0\}
\|_{L^2_{y_2}}(q_2)\right|\\
\lesssim
\int_{|q_1|\le M}  dq_1 ~ k^\chi(p,q_1)~
\left(\left(
\int_{\mathbb{R}^3} dq_2~
\left|  w_{j}(q_2)   k^\chi(q_1,q_2) \right|^2
\right)^{1/2}\right.\\
\left. \times
\left(
\int_{\mathbb{R}^3} dq_2~
\left|w_{-j}(q_2) \|\{  U(s_2) f_0 \}\|_{L^2_{y_2}}(q_2)\right|^2
\right)^{1/2} \right)  \\
\le
C_M \|w_{-j}f\|_{L^2_{q,x}}(s_2).
\end{multline*}
These estimates hold $\forall j > 0$.
Above we used the splitting as in \eqref{comparablePQ}.
Then furthermore, as in the previous estimates including \eqref{MdecayEST}, we have
\begin{multline*}
\|
w_\ell
H_5^{low,2}
\|_{L^\infty_pL^2_x}
\le
C_{M} 
\int_{\kappa}^t  ~ds_1 ~
\int_0^{s_1-\kappa} ~ ds_2 ~
 e^{-C(t- s_1)/M^{b/2} }
e^{-C(s_1- s_2)/M^{b/2} }
\\
\times \| w_{-j}\left\{  U(s_2) f_0 \right\}\|_{L^2_{p,x}}
\\
\le
C_{ M} 
\int_{0}^t  ~ds_1~   e^{-\frac{C}{2}(t- s_1)/M^{b/2} }
\int_0^{t} ~ ds_2 ~
e^{-\frac{C}{2}(t- s_1)/M^{b/2} }
e^{-\frac{C}{2}(s_1- s_2)/M^{b/2} }
\\
\times  \| w_{-j}\left\{  U(s_2) f_0 \right\}\|_{L^2_{p,x}}.
\end{multline*}
Notice that the first exponential controls the $s_1$ time integral, and the second and third exponential control the remaining time integral as follows
\begin{multline*}
\|
w_\ell
H_5^{low,2}
\|_{L^\infty_pL^2_x}\le
C_{ M} 
\int_0^{t} ~ ds_2 ~
e^{-\frac{C}{2}(t- s_2)/M^{b/2} }  \| w_{-j}\left\{  U(s_2) f_0 \right\}\|_{L^2_{p,x}}
 \\
\le
C_{ M} 
\int_0^{t} ~ ds_2 ~
e^{-\frac{C}{4}(t- s_2)/M^{b/2} }
  \| w_{-j}\left\{  U(s_2) f_0 \right\}\|_{L^2_{p,x}}.
\end{multline*}  
This completes the proof of Lemma \ref{hESTIMATE3}, by first choosing $M$ large. \qed\\

\section{Linear decay theory in $L^\infty_{p}L^\infty_x$}\label{sec:IIdecay}

Now we work on the linear $L^\infty_{p}L^\infty_x$ bounds and decay. The argument that we will use in this short section is not new. 
By Lemma 4.2, Lemma 4.3 and Lemma 4.4 from \cite{Strain2010}, for $\ell \geq 0, k\geq 0$, small $\eta$ and (possibly large) $j>0$, we have
\begin{multline*}
\|w_\ell\{U(t)f_0\} \|_{L^{\infty}_{p}L^\infty_x}
\\
\le C (1+t)^{-k}\| w_{\ell+k}f_0\|_{L^{\infty}_pL^\infty_x}+\eta(1+t)^{-k} ||\varpi_k w_\ell f ||_{L^\infty_{p,t}L^\infty_x}\\
+
C_\eta 
\int_0^t ds ~ e^{-\eta (t-s)} \|w_{-j}f\|_{L^2_{x,p}}(s).
\end{multline*}
Note that even though the lemmas just cited from \cite{Strain2010} were proven in the context of $\mathbb{T}^3_x$, they generalize directly to $\mathbb{R}^3_x$ without modification.

Using the same upper bound for $\int_0^t ds ~ e^{-\eta (t-s)} \|w_{-j}f\|_{L^2_{x,p}}(s)$ as in Lemma \ref{hESTIMATE3}, we then have the following result.

\begin{theorem}
\label{decay02}
Given $\ell \geq 0,\, r\in[1,2]$ and $k\in[0,\sigma_{r,0}]$.  Suppose that initially we have 
$
w_{\ell+k}f_0 \in L^\infty_{p}\left(L^2_x \cap L^\infty_x \right),
$
then under \eqref{hypSOFT} the semi-group satisfies
$$
\|w_\ell\{U(t)f_0\} \|_{L^{\infty}_{p}L^\infty_x}\\
\le
C (1+t)^{-k}
\left( \| w_{\ell+k}f_0\|_{L^\infty_{p}\left(L^2_x \cap L^\infty_x \right)}
+\|f_0\|_{L^2_pL^r_x} \right).
$$
Above the positive constant $C= C_{\ell, k}$ only depends on $\ell$ and $k$.  
\end{theorem}

\section{Nonlinear decay theory and global existence}\label{sec:NLdecay}
Suppose $f=f(t,x,p)$ solves \eqref{rBoltz0}
with initial condition $f(0,x,p) =  f_0(x,p)$.  
We may express mild solutions to this problem \eqref{rBoltz0} in the form
\begin{equation}
f(t,x,p) = \{U(t)f_0\}(x,p)+N[ f, f](t,x,p),
\label{nonLinProbSol}
\end{equation}
where we have used the notation
$$
N[ f_1, f_2](t,x,p)
\eqdef
\int_0^t ~ ds ~ \{U(t-s) \Gamma [ f_1(s), f_2(s)]\}(x,p).
$$
Here as usual $U(t)$ is the semi-group \eqref{Udef} which represents mild solutions to the linear problem 
\eqref{rBoltz}.  The purpose of this section is to prove Theorem \ref{mainGlobal}.

As a first step, we will use the following non-linear estimate.

\begin{lemma}
\label{nonlin0}  
Considering the non-linear operator defined in \eqref{gamma0} with \eqref{hypSOFT},
we have the following pointwise estimates  
$$
\| w_\ell \Gamma (h_1, h_2) \|_{L^1_x\cap L^\infty_x}(p)
\lesssim  \nu(p)
\| w_\ell h_1 \|_{L^\infty_q (L^2_x\cap L^\infty_x)}  \| w_\ell h_2 \|_{L^\infty_q (L^2_x\cap L^\infty_x)}. 
$$
These hold for any $\ell \ge 0$.
Furthermore,
$$
\| w_{\ell+1} \Gamma (h_1, h_2) \|_{L^\infty_p (L^1_x\cap L^\infty_x)}
\lesssim  
\| w_\ell h_1 \|_{L^\infty_p (L^2_x\cap L^\infty_x)}  \| w_\ell h_2 \|_{L^\infty_p (L^2_x\cap L^\infty_x)}. 
$$
\end{lemma}
Note that just as in the previous section, we use $q$ to denote the momentum variable in $\| w_\ell h_1 \|_{L^\infty_q (L^2_x\cap L^\infty_x)}  \| w_\ell h_2 \|_{L^\infty_q (L^2_x\cap L^\infty_x)}$ to avoid possible confusion, since $p$ is a variable in this inequality. 
The lemma above combined with Proposition \ref{BasicDecay} will be important tools in our proof of Theorem \ref{mainGlobal}. \\

\noindent {\it Proof of Lemma \ref{nonlin0}.}   We recall \eqref{nuDEF}, \eqref{gamma0}, and \eqref{weight}. 
For $\ell \ge 0$, it follows from \eqref{collisionalCONSERVATION} that
$$
w_\ell (p) 
\lesssim (p^0)^{\ell b/2}
\lesssim (p^{\prime 0})^{\ell b/2}  (q^{\prime 0})^{\ell b/2}
\lesssim w_\ell (p') w_\ell (q').
$$
A proof of this estimate above was given in \cite[Lemma 2.2]{MR1211782}.
Thus 
\begin{multline*}
w_\ell(p) \| \Gamma (h_1,h_2) \|_{L^1_x\cap L^\infty_x}\\
\lesssim
\int_{\mathbb{R}^3\times \mathbb{S}^{2}} ~ d\omega dq 
~ v_{\o} ~ \sigma(g,\theta)~
 \sqrt{J(q)} ~ 
 w_\ell (p') w_\ell (q') \| h_1\|_{L^2_x\cap L^\infty_x}(p^{\prime })\|h_2 \|_{L^2_x\cap L^\infty_x}(q^{\prime})
 \\
 +
\int_{\mathbb{R}^3\times \mathbb{S}^{2}} ~ d\omega dq 
~ v_{\o} ~ \sigma(g,\theta)~
 \sqrt{J(q)} ~ 
w_\ell(p)  \| h_1\|_{L^2_x\cap L^\infty_x}(p)\|h_2 \|_{L^2_x\cap L^\infty_x}(q)
 \\
\lesssim
  \| w_\ell h_1 \|_{L^\infty_q (L^2_x\cap L^\infty_x)}  \| w_\ell h_2 \|_{L^\infty_q (L^2_x\cap L^\infty_x)}
  \int_{\mathbb{R}^3\times \mathbb{S}^{2}} ~ d\omega dq 
~ v_{\o} ~ \sigma(g,\theta)~
 J^{1/2}(q)
  \\
\lesssim
 \nu(p)\| w_\ell h_1 \|_{L^\infty_q (L^2_x\cap L^\infty_x)}  \| w_\ell h_2 \|_{L^\infty_q (L^2_x\cap L^\infty_x)}.
\end{multline*}
The last inequality above follows directly from Lemma \ref{nuEST} since both the integral and $\nu(p)$ have the same asymptotic behavior at infinity.
That yields the first estimate.  
For the second estimate we notice 
from the first estimate that
$$
\| w_{\ell+1} \Gamma (h_1, h_2) \|_{L^1_x\cap L^\infty_x}(p)
\lesssim  w_1(p)\nu(p)
\| w_\ell h_1 \|_{L^\infty_q (L^2_x\cap L^\infty_x)}  \| w_\ell h_2 \|_{L^\infty_q (L^2_x\cap L^\infty_x)}.
$$
But $w_{1}(p)  \nu(p) \lesssim 1$ from 
Lemma \ref{nuEST} and \eqref{weight}.  This
 completes the proof.
\qed \\

\noindent {\it Proof of Theorem \ref{mainGlobal}}.  We will prove 
Theorem \ref{mainGlobal}
in three steps.   The first step gives existence, uniqueness and time decay via the contraction mapping argument.  The second step will establish continuity, and the last step shows positivity.

{\bf Step 1. Existence and Uniqueness}.  
When proving existence of mild solutions to \eqref{nonLinProbSol} it is natural to consider the mapping
$$
M[f]
\eqdef
\{U(t)f_0\}(x,p)+N[ f, f](t,x,p).
$$
We will show that this is a contraction mapping on  the space 
$$
M_{k,\ell}^{R} \eqdef \{f: [0,\infty) \times \mathbb{R}^3_x \times \mathbb{R}^3_p  \to \mathbb{R}\, |\,  
\|\varpi_k w_\ell f \|_{L^\infty_{p,t} (L^2_x\cap L^\infty_x)} \le R \},
\quad 
R>0.
$$
Here we recall the temporal weight \eqref{TIMEweight}.
We first estimate the non-linear term $N[ f, f]$ defined in the equation display below \eqref{nonLinProbSol}.
We apply Theorem \ref{decay0l} and Theorem \ref{decay02}, with $\ell \geq 0$,  and $k\in(1/2,\sigma_{r,0}]$  to obtain
\begin{multline*}
w_{\ell}(p)
\|
N[ f_1, f_2]
\|_{L^2_x\cap L^\infty_x}(t,p)\\
\lesssim
\int_0^t ~ ds ~
w_{\ell}(p)
\|
 \{U(t-s) \Gamma [ f_1(s), f_2(s)]\}
\|_{L^2_x\cap L^\infty_x}(p)
\\
\lesssim
\int_0^t ~  \frac{ds}{(1+t-s)^k}\left(
\left\|w_{\ell+k}
  \Gamma [ f_1(s), f_2(s)]
\right\|_{L^\infty_q({L^2_x\cap L^\infty_x})}+\left\|  \Gamma [ f_1(s), f_2(s)]
\right\|_{L^2_qL^r_x}\right).
\end{multline*}
When $\ell>3/b-1$ (since $k\le 1$), using interpolation the above is bounded by  
$$\int_0^t ~  \frac{ds}{(1+t-s)^k}
\left\|w_{\ell+1}
  \Gamma [ f_1(s), f_2(s)]
\right\|_{L^\infty_q({L^1_x\cap L^\infty_x})}.$$
Next Lemma \ref{nonlin0} allows us to bound the above by
\begin{gather*}
\lesssim
\int_0^t ~  \frac{ds}{(1+t-s)^k}
\left\| w_\ell f_1(s) \right\|_{L^\infty_q({L^2_x\cap L^\infty_x})}
\left\| w_\ell  f_2(s) \right\|_{L^\infty_q({L^2_x\cap L^\infty_x})}.
\end{gather*}
From Proposition \ref{BasicDecay} and \eqref{TIMEweight} we see that the last line is
\begin{gather*}
\lesssim
|| \varpi_k w_\ell f_1 ||_{L^\infty_{p,t}({L^2_x\cap L^\infty_x})}
|| \varpi_k w_\ell f_2 ||_{L^\infty_{p,t}({L^2_x\cap L^\infty_x})}
\int_0^t ~  \frac{ds}{(1+t-s)^k (1+s)^{2k}}
\\
\lesssim
(1+t)^{-k}
|| \varpi_k w_\ell f_1 ||_{L^\infty_{p,t}({L^2_x\cap L^\infty_x})}
|| \varpi_k w_\ell f_2 ||_{L^\infty_{p,t}({L^2_x\cap L^\infty_x})}.
\end{gather*}
We have shown 
$$
||  \varpi_kw_\ell N[ f_1, f_2] ||_{L^\infty_{p,t}({L^2_x\cap L^\infty_x})}
\lesssim
|| \varpi_k w_\ell f_1 ||_{L^\infty_{p,t}({L^2_x\cap L^\infty_x})}
|| \varpi_k w_\ell f_2 ||_{L^\infty_{p,t}({L^2_x\cap L^\infty_x})}.
$$
To handle the linear semigroup,  $U(t)$, we again use Theorem \ref{decay0l} and Theorem \ref{decay02} to obtain
\begin{multline*}
|| \varpi_k w_\ell M[ f] ||_{L^\infty_{p,t}({L^2_x\cap L^\infty_x})}\\
\le
 C_{\ell,k}\left(\left\| w_{\ell+k} f_0 \right\|_{L^\infty_{p}({L^2_x\cap L^\infty_x})}+\left\|  f_0 \right\|_{L^2_{p}L^r_x}
 +
|| \varpi_k w_\ell f ||_{L^\infty_{p,t}({L^2_x\cap L^\infty_x})}^2
 \right).
\end{multline*}
We conclude that $M[\cdot]$ maps $M_{k,\ell}^{R}$ into itself for $0<R$ chosen sufficiently small and e.g.
$
\left\| w_{\ell+k} f_0 \right\|_{L^\infty_{p}({L^2_x\cap L^\infty_x})}+\left\|  f_0 \right\|_{L^2_{p}L^r_x}\le \frac{R}{2 C_{\ell,k}}.
$
To obtain a contraction, we consider the difference
$$
M[ f_1] - M[ f_2] 
=
N[ f_1-f_2, f_1]
+
N[ f_2, f_1-f_2].
$$
Then as in the previous estimates we have
\begin{multline*}
||  \varpi_k w_\ell (M[ f_1] - M[ f_2])||_{L^\infty_{p,t}({L^2_x\cap L^\infty_x})} \\
\le
C_{\ell,k}^*
\left( 
||  \varpi_k w_\ell f_1 ||_{L^\infty_{p,t}({L^2_x\cap L^\infty_x})}
 +
|| \varpi_k w_\ell  f_2 ||_{L^\infty_{p,t}({L^2_x\cap L^\infty_x})}
 \right)\\
\times || \varpi_k w_\ell( f_1 - f_2) ||_{L^\infty_{p,t}({L^2_x\cap L^\infty_x})}.
\end{multline*}
With these estimates, the existence and uniqueness of solutions to \eqref{rBoltz} follows from the contraction mapping principle on $M_{k,\ell}^{R}$ when $R>0$ is suitably small.

 {\bf Step 2. Continuity}.  We perform the estimates from Step 1 on the space
 $$
M_{k,\ell}^{R,0} \eqdef  C^0 ([0,\infty) \times \mathbb{R}^3_x \times \mathbb{R}^3_p) \cap M_{k,\ell}^{R},
\quad 
R>0.
$$
As in Step 1, 
we have a uniform in time contraction mapping on $M_{k,\ell}^{R,0}$ for suitable $R$.  
Furthermore $M[f]$ is continuous if $f\in M_{k,\ell}^{R,0}$  and $f_0$ is continuous.  
Since the convergence is uniform, the limit will be continuous  globally in time.
This argument is standard and we refer for instance to \cite{MR1379589,MR1211782,MR2679358} for full details.

{\bf Step 3. Positivity}. 
We use the standard alternative approximating formula
$$
\left(\partial_t +\hat{p}\cdot \nabla_x \right) F^{n+1} + R(F^{n})F^{n+1}
=
\mathcal{Q}_+(F^{n},F^{n}).
$$
with the same initial conditions 
$
\left. F^{n+1}\right|_{t=0} =  F_0= J + \sqrt{J} f_0,
$
for $n\ge 1$ and for instance
$
 F^{1} \eqdef J + \sqrt{J} f_0.
$
Here we have used the standard decomposition of the collision operator
$
\mathcal{Q}
=
\mathcal{Q}_+ - \mathcal{Q}_-
$
into gain and loss terms with
$$
\mathcal{Q}_-(F^{n+1},F^{n}) 
=
\int_{\mathbb{R}^3}  dq ~
\int_{\mathbb{S}^{2}} d\omega ~
~ v_{\o} ~ \sigma (g,\theta ) ~ F^{n+1}(p) F^{n}(q)
= 
R(F^{n})F^{n+1},
$$
and
$
R(F^{n}) \eqdef \mathcal{Q}_-(1,F^{n}).
$
If we consider 
$
 F^{n+1}(t,x,p) = J + \sqrt{J} f^{n+1}(t,x,p),
$
then related to Step 1 we may show that 
$
w_\ell(p)f^{n+1}(t,x,p)
$
is convergent in $L^\infty_{p,t}(L^2_x\cap L^\infty_x)$ on a local time interval $[0,T]$ where $T$ will generally depend upon the size of the initial data. 
In particular $f^{n+1}(t,x,p) = \frac{ F^{n+1} - J}{\sqrt{J}}$ satisfies the equation
$$
\left(\partial_t +\hat{p}\cdot \nabla_x  + \nu(p)\right)f^{n+1}
= K(f^{n})
+
\Gamma_+(f^{n},f^{n})
-
\Gamma_- (f^{n+1},f^{n}).
$$
Here again we use the standard definition
$
\Gamma_+
-
\Gamma_- = \Gamma
$ 
where
$$
\Gamma_-(f^{n+1},f^{n})=
\int_{\mathbb{R}^3}  dq ~
\int_{\mathbb{S}^{2}} d\omega ~
~ v_{\o} ~ \sigma (g,\theta )  \sqrt{J(q)} ~ f^{n+1}(p) f^{n}(q).
$$
We rewrite this equation using the solution formula to the system \eqref{rBlinWO} as
 $$
f^{n+1}
= 
G(t) f_0
+
\mathcal{L}(f^{n+1}, f^{n}).
$$
This solution formula $G(t)$ is defined just below \eqref{rBlinWO}.  Furthermore
\begin{multline*}
\mathcal{L}(f^{n+1}, f^{n})
\eqdef
 \int_0^t ~ ds ~ 
G(t-s) K(f^{n}) 
\\
+ \int_0^t ~ ds ~ 
G(t-s)  \Gamma_+(f^{n},f^{n})  
-
G(t-s)  \Gamma_- (f^{n+1},f^{n}).
\end{multline*}
For given $T>0$ and  $R > 0$ we consider the space $M_{k,\ell}^{R}([0,T])$ defined by
 $$
 \left\{f : [0,T] \times \mathbb{R}^3_x \times \mathbb{R}^3_p\to\mathbb{R}\,|\, \| \varpi_k w_\ell f \|_{L^\infty_{[0,T]}L^\infty_p (L^2_x\cap L^\infty_x)}\le R \right\}.
$$
We can do something precisely similar to Step 1 to prove the existence  of  
$f^{n+1} \in  M_{k,\ell}^{R}([0,T])$. Instead of using Theorem \ref{decay0l} and Theorem \ref{decay02} to deal with the operator $U$, we use in particular Lemma \ref{hESTIMATES} in this paper and Lemma 4.2 in \cite{Strain2010} to deal with the operator $G$. Then we have estimates similar to those in Step 1 for our mapping in the space $M_{k,\ell}^{R}([0,T])$.
Now given 
$f^{n} \in  M_{k,\ell}^{R}([0,T])$
and
$
\left\|  w_{\ell+k}f_0 \right\|_{L^\infty_p (L^2_x\cap L^\infty_x)}+\left\| f_0 \right\|_{L^2_p L^r_x} \le \frac{R}{2 C_{k,\ell}}
$
with $R>0$ chosen sufficiently small, as in Step 1, we can prove the existence  of  
$f^{n+1} \in  M_{k,\ell}^{R}([0,T])$.

With the estimates established in this paper, it is now not hard to show that
$$
\| w_\ell (f^{n+1} - f^{n}) \|_{L^\infty_{[0,T]}L^\infty_p (L^2_x\cap L^\infty_x) } \le C T  \| w_\ell (f^{n} - f^{n-1}) \|_{L^\infty_{[0,T]}L^\infty_p (L^2_x\cap L^\infty_x) }).
$$ 
Here
 $T>0$  is sufficiently small, and
 the constant $C>0$ can be chosen independent of any small $T$.  Therefore there exists a $T^*>0$ such that  $w_\ell f^{n} \to w_\ell f$ uniformly in $L^\infty_p (L^2_x\cap L^\infty_x)$ on $[0,T^*]$. 
This will be sufficient to prove the positivity globally in time.  

Indeed if $F^{n} \ge 0$, then so is $\mathcal{Q}_+(F^{n},F^{n})\ge 0$.  With the  representation formula
\begin{multline*}
F^{n+1}(t,x,p) 
=
e^{-\int_0^t  ds  R(F^{n})(s, x- \hat{p}(t-s), p)} ~
F_0(x-\hat{p}t, p)
\\
+
\int_0^t  ds ~ e^{-\int_s^t  d\tau  R(F^{n})(\tau, x- \hat{p}(t-\tau), p)} ~ \mathcal{Q}_+(F^{n},F^{n})(s, x- \hat{p}(t-s), p).
\end{multline*}
Induction shows  $F^{n+1}(t,x,p)  \ge 0$ for all $n\ge 0$ if $F_0\ge 0$,
which implies in the limit $n\to\infty$ that $F(t,x,p)= J + \sqrt{J} f(t,x,p)\ge 0$.  Using our uniqueness, this is the same  $F$ as the one from Step 1 on the time interval $[0,T^*]$.  We  extend this positivity for all time intervals $[0,T^*] + T^* k$ for any $k\ge 1$ by repeating this  procedure and using the global uniform bound in $M_{k,\ell}^{R}$ from Step 1.
\qed

\end{document}